\documentclass[10pt,a4paper]{elsarticle}

\usepackage{float}
\usepackage{amstext,amsfonts,amsmath,amssymb}
\usepackage{amsthm}
\usepackage{amsmath}
\usepackage{breqn}
\usepackage{makeidx}
\usepackage{bm}

\makeatletter
\def\ps@pprintTitle{%
  \let\@oddhead\@empty
  \let\@evenhead\@empty
  \def\@oddfoot{\reset@font\hfil\thepage\hfil}
  \let\@evenfoot\@oddfoot
}
\def\@author#1{\g@addto@macro\elsauthors{\normalsize%
\def\baselinestretch{1}%
\upshape\authorsep#1\unskip\textsuperscript{%
\ifx\@fnmark\@empty\else\unskip\sep\@fnmark\let\sep=,\fi
\ifx\@corref\@empty\else\unskip\sep\@corref\let\sep=,\fi}%
\def\authorsep{\unskip,\space}%
\global\let\@fnmark\@empty
\global\let\@corref\@empty 
\global\let\sep\@empty}%
\@eadauthor={#1}
}
\makeatother

\newtheorem{theorem}{Theorem}[section]

\newtheorem{lemma}[theorem]{Lemma}

\def\nfrac#1#2{\mbox{\footnotesize $\displaystyle\frac{#1}{#2}$}}

\begin{document}

\centerline {\bf \large ACCURATE APPROXIMATIONS OF SOME EXPRESSIONS}

\smallskip

\centerline {\bf \large INVOLVING TRIGONOMETRIC FUNCTIONS}

\begin{center}
Marija Nenezi\' c${}^{\,\mbox{\tiny 1)}}$,
Branko Male\v sevi\' c${}^{\mbox{\tiny 1)}}$,
Cristinel Mortici${}^{\mbox{\tiny *},\mbox{\tiny 2)}}$
\end{center}

\footnotetext{$\!\!\!\!\!\!\!\!\!\!\!{}^{*}\,$Corresponding author
}

\footnotetext{$\!\!\!\!\!\!\!\!\!\!\!$E-mails:
Marija~Nenezi\' c$\,<${\sl maria.nenezic@gmail.com}$>$, Branko~Male\v sevi\' c$\,<${\sl malesevic@etf.rs}$>$,
Cristinel~Mortici$\,<${\sl cmortici@valahia.ro}$>$}

\begin{center}
{\footnotesize \it
${}^{1)}$Faculty of Electrical Engineering, University of Belgrade,                           \\[+0.20 ex]
Bulevar kralja Aleksandra 73, 11000 Belgrade, Serbia                                          \\[+1.00 ex]
${}^{2)}$Valahia University of T\^{a}rgovi\c{s}te, Bd. Unirii 18, 130082 T\^{a}rgovi\c{s}te;  \\[-0.80 ex]
Academy of Romanian Scientists, Splaiul Independen \c{t}ei 54, 050094 Bucharest, Romania}
\end{center}

\bigskip
\noindent
{\small {\bf Abstract.}
The aim of this paper is to apply an original computation method due to Male\v sevi\' c and Makragi\' c \cite{[7]} to the problem of approximating some trigonometric functions. Inequalities of Wilker-Cusa-Huygens are discussed, but the method can be successfully applied to a wide class of problems. In particular, we improve the estimates recently obtained by Mortici \cite{[1]} and moreover we show that they hold true also on some extended intervals.}

\bigskip
\noindent
{\footnotesize {\bf Keyword.}
Wilker inequality, trigonometric approximation}

\section{Introduction} 
\smallskip
\noindent
In the reference \cite{[2]} {\sc J.$\,$B. Wilker} presented the inequality
\begin{equation}
\label{Ineq_1}
2<\left(\nfrac{\sin x}{x}\right)^{\!2} + \nfrac{\tan x}{x},
\end{equation}
for $x \!\in\! \left(0, \pi/2\right)$ and he asked for largest constant $c$ in
\begin{equation}
\label{Ineq_2}
2+c x^3 \tan x
<
\left(\nfrac{\sin x}{x} \right)^{\!2} + \nfrac{\tan x}{x}\quad \mbox {for} \quad c>0,
\end{equation}
and for $x \!\in\! \left(0, \pi/2\right)$. Recently, Wilker inequality is a lot studied
in different paper works. In the paper \cite{[4]}, {\sc J.S.Sumner}, {\sc A.A.Jagers}, {\sc M. Vowe},
{\sc J. Anglesio} proved the following double inequality
\begin{equation}
\label{Ineq_3}
2+\nfrac{16}{\pi^4} x^3 \tan x
<
\left(\nfrac{\sin x}{x} \right)^{\!2} + \nfrac{\tan x}{x}
<
2+\nfrac{8}{45} x^3 \tan x,
\end{equation}
for $x \!\in\! \left(0, \pi/2\right)$. In the paper \cite{[1]}, {\sc C. Mortici} has proved
the~following~two statements:
\begin{theorem}
For every $x \!\in\! \left(0, 1\right)$ we have$:$
\begin{equation}
\label{Ineq_4}
2+\left(\nfrac{8}{45}-a\left(x \right)\right) x^3 \tan x
<
\left(\nfrac{\sin x}{x}\right)^{\!2} + \nfrac{\tan x}{x}
<
2+\left(\nfrac{8}{45}-b\left(x \right)\right) x^3 \tan x,
\end{equation}
where $a\left(x \right)= \nfrac{8}{945} x^2$,
$ b\left(x \right)= \nfrac{8}{945} x^2- \nfrac{16}{14175} x^4$.
\end{theorem}
\begin{theorem}
For every $x \!\in\! \left(\nfrac{\pi}{2}-\nfrac{1}{2}, \nfrac{\pi}{2} \right)$
in the left-hand side and for every $x \!\in\! \left(\nfrac{\pi}{3}-\nfrac{1}{2}, \nfrac{\pi}{2} \right)$
in the right-hand side the following inequalities are true$:$
\begin{equation}
\label{Ineq_5}
2+\left(\nfrac{16}{\pi^4}+c\left(x \right)\right)x^3 \tan x
<
\left(\nfrac{\sin x}{x} \right)^{\!2}+\nfrac{\tan x}{x}
<
2+\left(\nfrac{16}{\pi^4}+d\left(x \right) \right)x^3 \tan x ,
\end{equation}
where \\

\vspace*{-2.5 mm}

\noindent
$$
c\left(x\right)
\!=\!
\left(\nfrac{160}{\pi^5}\!-\!\nfrac{16}{\pi^3}\right)\!\left(\nfrac{\pi}{2}\!-\!x\right)
,
~
d\left(x\right)
\!=\!
\left(\nfrac{160}{\pi^5} \!-\!\nfrac{16}{\pi^3}\right)\!\left(\nfrac{\pi}{2}\!-\!x\right)\!+\!\left(\nfrac{960}{\pi^6}\!-\!\nfrac{96}{\pi^4}\right)\!\left(\nfrac{\pi}{2}\!-\!x\right)^{\!2}.
$$
\end{theorem}

\noindent

\noindent
Theorem 1.1. and Theorem 1.2. describe a subtly analysis of Wilker inequality by {\sc C. Mortici}.
The method of proving inequalites in this paper was given in the paper \cite{[7]}
and it is based on use of appropriate approximations of some mixed trigonometric polynomials with
finite Taylor series. The method presents continuation of method of {\sc C. Mortici} presented in \cite{[6]}.
The method from the paper \cite{[7]} was applied in papers \cite{[8]} and \cite{[9]} on inequalities which are
closely related.

\section{The Main Results} 
\smallskip
\noindent
The main purpose of our paper is to extend the intervals defined in theorems given by {\sc C. Mortici}
\cite{[6]}. More precisely, we extend the domains $\left(0, 1\right)$ and $\left(\nfrac{\pi}{2}-\nfrac{1}{2}, \nfrac{\pi}{2} \right)$ from the previous theorems to $\left(0, \nfrac{\pi}{2} \right)$. We give the next two statements.
\begin{theorem}
For every $x \!\in\! \left(0, \nfrac{\pi}{2} \right)$ the following inequalities are true$:$
\begin{equation}
\label{Ineq_7}
2+\left(\nfrac{8}{45}-a\left(x \right) \right) x^3 \tan x
<
\left(\nfrac{\sin ⁡x}{x} \right)^{\!2}+\nfrac{\tan x}{x}
<
2+\left(\nfrac{8}{45}-b_{1} \left(x \right)\right) x^3 \tan x,
\end{equation}
where $a(x)= \nfrac{8}{945} x^2$, $b_{1}(x)= \nfrac{8}{945} x^2- \nfrac{\mbox{\boldmath $a$} }{14175} x^4$
and
$
\mbox{\boldmath $a$}
=
\frac{480\pi^6-40320\pi^4+3628800}{\pi^8}
=
17.15041 \ldots \, .
$
\end{theorem}
\begin{theorem}
For every $x \!\in\! \left(0, \nfrac{\pi}{2} \right)$ the following inequalities are true$:$
\begin{equation}
\label{Ineq_8}
2+\left(\nfrac{16}{\pi^4} +c\left(x\right) \right) x^3 \tan x
<
\left(\nfrac{\sin x}{x}\right)^{\!2}+\nfrac{\tan x}{ x}
<
2+\left(\nfrac{16}{\pi^4}+d\left(x\right) \right) x^3 \tan x,
\end{equation}
where \\

\vspace*{-3.0 mm}

\noindent
$$
c(x)
\!=\!
\left(\nfrac{160}{\pi^5}\!-\!\nfrac{16}{\pi^3}\right)\!\left(\nfrac{\pi}{2}\!-\!x\right),
~
d(x)
\!=\!
\left(\nfrac{160}{\pi^5}\!-\!\nfrac{16}{\pi^3}\right)\!\left(\nfrac{\pi}{2}\!-\!x\right)+\left(\nfrac{960}{\pi^6}-\nfrac{96}{\pi^4}\right)\!\left(\nfrac{\pi}{2}\!-\!x\right)^{\!2}.
$$
\end{theorem}
\noindent
In \cite{[7]} is considered a method of proving trigonometric inequalities for mixed trigonometric polynomials:
\begin{equation}
\label{Ineq_9}
f(x)=\sum\limits_{i=1}^{n}\alpha_{i}x^{p_{i}}\!\cos^{q_{i}}\! x\sin^{r_{i}}\! x>0,
\end{equation}
for $x \!\in\!(\delta_{2},\delta_{1})$, $\delta_{2}\!<\!0\!<\!\delta_{1}$, where $\alpha_{i}
\!\in\! \mathbb{R} \!\setminus\! \lbrace0\rbrace$, \mbox{$p_{i},\;q_{i},\;r_{i} \!\in\! \mathbb{N}_{0}$}
and \mbox{$n \!\in\! \mathbb{N}$}. One method of proving inequalities in form (\ref{Ineq_9})
is based on transformation, using the sum of sine and cosine of multiple angles.

\smallskip
\noindent
Let us mention some facts from~\cite{[7]}. Let $\varphi : [a,b] \longrightarrow \mathbb{R}$
be a function which is differentiable on a segment $[a,b]$ and differentiable arbitrary
number of times on a right neighbourhood of the point $x=a$ and denote by $T^{\varphi, a}_{m}(x)$
the Taylor polynomial of the function $\varphi(x)$ in the point $x=a$ of the order $m$. If~there
is some $\eta\!>\!0$ such that holds:
\mbox{$
T^{\varphi, a}_{m}(x) \geq \varphi(x)
$}, for $x \!\in\! (a,a+\eta) \!\subset\! [a,b]$;
then let us define $\overline{T}^{\,\varphi,a}_{m}(x)=T^{\,\varphi, a}_{m}(x)$ and
$\overline{T}^{\,\varphi,a}_{m}(x)$ present an upward approximation of the function $\varphi(x)$ on right
neighbourhood $(a,a+\eta)$ of the point $a$ of the oreder~$m$. Analogously, if there is some $\eta>0$
such that holds:
\mbox{$
T^{\varphi, a}_{m}(x) \leq \varphi(x)
$}, for $x \!\in\! (a,a+\eta) \!\subset\! [a,b]$;
then let us define $\underline{T}^{\varphi,a}_{\,m}(x)=T^{\,\varphi, a}_{m}(x)$
and $\underline{T}^{\varphi,a}_{\,m}(x)$ present a downward approximation of the function $\varphi(x)$
on right neighbourhood $(a,a+\eta)$ of the point $a$ of the order $m$. Let us note that it is possible
to analogously define upward and downward approximations on some left neighbourhood of a point.

\noindent
According to the paper \cite{[7]} following Lemmas are true:
\begin{lemma}
{\rm \textbf{\!(i)}} For the polynomial $T_n(t)=\!\!\!\!\displaystyle\sum\limits_{i=0}^{(n-1)/2}
\!\!\dfrac{(-1)^it^{2i+1}}{(2i+1)!}$, where $n=4k+1$,\\

\vspace*{-2.5 mm}

\noindent
$k\in \mathbb{N}_{0}$, it is valid:
\begin{equation}
\label{Ineq_10}
\Big{(}\forall t \in \big{[}0, \mbox{\small $\sqrt{(n+3)(n+4)}$}\,\big{]}\Big{)}\,\overline{T}_n(t)
\geq
\overline{T}_{n+4}(t)\geq \sin t,
\end{equation}
\begin{equation}
\label{Ineq_11}
\Big{(}\forall t \in \big{[}\mbox{\small $-\sqrt{(n+3)(n+4)}$},0\big{]}\Big{)}\,\underline{T}_n(t)
\leq
\underline{T}_{n+4}(t)\leq \sin t.
\end{equation}

\smallskip
\noindent
For the value $t=0$ the inequalities in {\rm $($\ref{Ineq_10}$)$} and {\rm $($\ref{Ineq_11}$)$} turn into equalities.
For the values $t\!=\!\mbox{\small $\pm\sqrt{(n+3)(n+4)}$}$ the equalities
$\overline{T}_n(t)\!=\!\overline{T}_{n+4}(t)$ and $\underline{T}_n(t)\!=\!\underline{T}_{n+4}(t)$
are true, respectively.

\medskip
\noindent
{\rm \textbf{(ii)}}
For the polynomial $T_n(t)=\!\!\!\displaystyle\sum\limits_{i=0}^{(n-1)/2} \!\!\dfrac{(-1)^it^{2i+1}}{(2i+1)!}$,
where $n=4k+3$, $k\in \mathbb{N}_{0}$, it is \\

\vspace*{-2.5 mm}

\noindent
valid:

\vspace*{2.5 mm}

\noindent
\begin{equation}
\label{Ineq_12}
\Big{(}\forall t \in \big{[}0,\mbox{\small $\sqrt{(n+3)(n+4)}$}\,\big{]}\Big{)}\,\underline{T}_n(t)
\leq
\underline{T}_{n+4}(t)\leq \sin t,
\end{equation}
\begin{equation}
\label{Ineq_13}
\Big{(}\forall t \in \big{[}\mbox{\small $-\sqrt{(n+3)(n+4)}$},0\big{]}\Big{)}\,\overline{T}_n(t)
\geq
\overline{T}_{n+4}(t)\geq \sin t.
\end{equation}

\smallskip
\noindent
For the value $t=0$ the inequalities in {\rm $($\ref{Ineq_12}$)$} and {\rm $($\ref{Ineq_13}$)$} turn into equalities.
For the values $t=\mbox{\small $\pm\sqrt{(n+3)(n+4)}$}$ the equalities $\underline{T}_n(t)=\underline{T}_{n+4}(t)$
and $\overline{T}_n(t)=\overline{T}_{n+4}(t)$ are true, respectively.
\end{lemma}

\noindent
Let us notice that for the function $\sin x$ we have following order:
\begin{equation}
\label{Ineq_14}
\begin{array}{ll}
\underline{T}_{3}^{\sin ,0}(x) \leq \underline{T}_{7}^{\sin ,0}(x) \leq \underline{T}_{11}^{\sin ,0}(x) \leq \underline{T}_{15}^{\sin ,0}(x)
\,\leq\,\ldots\,\leq\,\sin x\,\leq
\ldots &                                                                                                 \\[1.5 ex]
\leq
\overline{T}_{13}^{\sin ,0}(x) \leq \overline{T}_{9}^{\sin ,0}(x) \leq \overline{T}_{5}^{\sin ,0}(x) \leq \overline{T}_{1}^{\sin ,0}(x)
\;\; \mbox{for} \;\; x\!\in\!\left[\,0,\sqrt{20}\,\right].
\end{array}
\end{equation}
\begin{lemma}
{\rm \textbf{(i)}}
For the polynomial $T_n(t)=\displaystyle\sum\limits_{i=0}^{n/2}\dfrac{(-1)^it^{2i}}{(2i)!}$, where $n\!=\!4k$, $k\!\in\!\mathbb{N}_{0}$, \\

\vspace*{-2.5 mm}

\noindent
it is valid:
\begin{equation}
\label{Ineq_15}
\Big{(}\forall t \in \big{[}\mbox{\small $-\sqrt{(n+3)(n+4)},\sqrt{(n+3)(n+4)}$}\,\big{]}\Big{)} \,\,\,\overline{T}_n(t)\geq \overline{T}_{n+4}(t)\geq \cos t.
\end{equation}

\smallskip
\noindent
For the value $t=0$ the inequality in {\rm $($\ref{Ineq_15}$)$} turns into equality.
For the values $t=\mbox{\small $\pm\sqrt{(n+3)(n+4)}$}$ the equality
$\overline{T}_n(t)=\overline{T}_{n+4}(t)$ is true.

\medskip
\noindent
{\rm \textbf{(ii)}}
For the polynomial $T_n(t)=\displaystyle\sum\limits_{i=0}^{n/2}\dfrac{(-1)^it^{2i}}{(2i)!}$,
where $n\!=\!4k\!+\!2$, $k\!\in\!\mathbb{N}_{0}$, it~is~valid:
\begin{equation}
\label{Ineq_16}
\Big{(}\forall t \in \big{[}\mbox{\small $-\sqrt{(n+3)(n+4)},\sqrt{(n+3)(n+4)}$}\,\big{]}\Big{)}
\,\,\, \underline{T}_n(t)\leq \underline{T}_{n+4}(t)\leq \cos t.
\end{equation}

\smallskip
\noindent
For the value $t=0$ the inequality in {\rm $($\ref{Ineq_16}$)$} turns into equality.
For the values $t=\mbox{\small $\pm\sqrt{(n+3)(n+4)}$}$ the equality $\underline{T}_n(t)=\underline{T}_{n+4}(t)$
is true.
\end{lemma}

\noindent
Let us notice that for the function $\cos x$ we have following order:
\begin{equation}
\label{Ineq_17}
\begin{array}{ll}
\underline{T}_{2}^{\cos ,0}\!(x)\!\leq \!\underline{T}_{6}^{\cos ,0}\!(x)\!\leq \!\underline{T}_{10}^{\cos ,0}\!(x)\!\leq \!\underline{T}_{14}^{\cos ,0}\!(x)
\,\leq\,\ldots\,\leq\,\cos x\,\leq\,\ldots
&  \\[1.5 ex]
\leq \overline{T}_{12}^{\cos ,0}(x)\!\leq\!\overline{T}_{8}^{\cos ,0}(x)\!\leq\!\overline{T}_{4}^{\cos ,0}(x)\!\leq\!\overline{T}_{0}^{\cos ,0}(x)\ \mbox{for} \ x\!\in\!\left[\,0,\sqrt{12}\,\right].
\end{array}
\end{equation}

\noindent
Proofs of previous Lemmas given above are presented in the paper \cite{[9]}.

\section{Proofs} 

In order to prove Theorem 2.1. and Theorem 2.2. we will separately observe left and right sides of inequalities.

\bigskip

\textbf{The proof of Theorem 2.1.}

\bigskip

{\small

Transforming inequality (\ref{Ineq_7}) we have following considerations.

\medskip
{\normalsize \bf(A)} Proving the left side of inequality
\begin{equation}
\label{Ineq_18}
2+\left(\nfrac{8}{45}-a(x) \right) x^3 \tan x<\left(\nfrac{\sin ⁡x}{x} \right)^{\!2}+\nfrac{\tan x}{x},
\end{equation}
for $x \!\in\! \left(0,\nfrac{\pi}{2}\right)$. The inequality (\ref{Ineq_18}) is equivalent to the
mixed trigonometric inequality
\begin{equation}
\label{Ineq_19}
\begin{array}{lcl}
f(x)
\!\!&\!\!=\!\!&\!\!
1\!-\!8 x^2\!+h_1(x)\cos 4 x + h_2(x)\cos2 x + h_3(x)\sin2 x \\[1.0 ex]
\!\!&\!\!=\!\!&\!\!
1\!-\!8 x^2\!-\!\cos4 x\!-\!8 x^2 \!\cos2 x\!+\!\left(4 x\!-\!4\left(\nfrac{8}{45}\!-\!a(x)\right) x^5 \right)\sin2 x\!>\!0,
\end{array} \end{equation}
for
$x \!\in\!\left(0,\nfrac{\pi}{2} \right)$,
and
$h_1(x)\!=\!-1<0$,
$h_2(x)\!=\!-8 x^2<0$,
$h_3(x)\!=\!4x-4\left(\nfrac{8}{45}-a(x)\right)x^5$.

\medskip
\noindent
Now let us consider two cases:

\medskip
\noindent
{\normalsize \bf (A/I)}
{\boldmath $x\!\in\!(0,1.57]$}
Let us determine sign of the polynomial $h_3(x)$. As we see, that polynomial is the polynomial of $7^{th}$ degree
\begin{equation}
\label{Ineq_21}
h_3(x)
=
P_7(x)
=
4 x-4\left(\nfrac{8}{45}-a(x) \right) x^5
=
\nfrac{32}{945} x^7 - \nfrac{32}{45} x^5 + 4 x.
\end{equation}
Using the factorization of the polynomial $P_7(x)$ we have
\begin{equation}
\label{Ineq_22}
P_7(x)=\nfrac{4}{945} x(8 x^6 -168 x^4+945)=\nfrac{4}{945} x P_6(x) , \end{equation}
where
\begin{equation}
\label{Ineq_23}
P_6(x)=8 x^6 -168 x^4+945,
\end{equation}
for $x\!\in\!(0,1.57]$. Introducing the substitution $s= x^2$ we can notice that the polynomial
$P_6(x)$ can be transformed into polynomial of $3^{rd}$ degree
\begin{equation}
\label{Ineq_24}
P_3(s)=8 s^3-168 s^2+945,
\end{equation}
for $s\!\in\!(0,2.4649]$. Using MATLAB software we can determine the real numerical factorization of the polynomial
\begin{equation}
\label{Ineq_25}
P_3(s)=\alpha(s - s_1)(s - s_2)(s - s_3), \end{equation}
where $\alpha = 8$ and where

\noindent
\begin{equation}
\label{Ineq_26}
\begin{array}{rcl}
s_1 \!\!&\!=\!&\!\!\!-\nfrac{1}{8}\left(18172\!+\!84 I \sqrt{21495}\right)^{1/3}\!-\!\nfrac{98}{\left(18172\!+\!84 I \sqrt{21495}\right)^{1/3}} \!+\! 7                                \\[-0.1 ex]
\!\!&\! \!&\!\!\!+\nfrac{3}{4} I \sqrt{3}\left(\nfrac{1}{6}\left(18172\!+\!84 I \sqrt{21495}\right)^{1/3} \!-\! \nfrac{392}{3\left(18172 \!+\! 84 I \sqrt{21495}\right)^{1/3}}\right)  \\[0.00 ex]
\!\!&\!=\!&\!\!-2.253\ldots ,                                                                                                                                                          \\[1.50 ex]
s_2 \!\!&\!=\!&\!\! \!-\nfrac{1}{8}\left(18172\!+\!84 I \sqrt{21495}\right)^{1/3}\!-\!\nfrac{98}{\left(18172\!+\!84 I \sqrt{21495}\right)^{1/3}} \!+\! 7                               \\[-0.1 ex]
\!\!&\! \!&\!\!\!-\nfrac{3}{4} I \sqrt{3}\left(\nfrac{1}{6}\left(18172\!+\!84 I \sqrt{21495}\right)^{1/3} \!-\! \nfrac{392}{3\left(18172 \!+\! 84 I \sqrt{21495}\right)^{1/3}}\right)  \\[0.00 ex]
\!\!&\!=\!&\!\! 2.528\ldots ,                                                                                                                                                          \\[1.50 ex]
s_3 \!\!&\!=\!&\!\!\nfrac{1}{4}\left(18172\!+\!84 I \sqrt{21495}\right)^{1/3} \!+\! \nfrac{196}{\left(18172 \!+\! 84 I \sqrt{21495}\right)^{1/3}} \!+\!7                               \\[0.00 ex]
\!\!&\!=\!&\!\! 20.724\ldots ;                                                                                                                                                         \\[1.50 ex]
\end{array}
\end{equation}
for $I= \sqrt{-1}$ (imaginary unit). The polynomial $P_3(s)$ has exactly three simple real roots with
a symbolic radical representation and corresponding numerical values $s_1$, $s_2$, $s_3$ given at (\ref{Ineq_26}).
Since $P_3(0)>0$ it follows that $P_3(s)>0$ for $s\!\in\! \left(s_1, s_2 \right)$, so we have following
conclusions:
\begin{equation}
\label{Ineq_27}
\begin{array}{rccl}
                & P_3(s)\!>\!0 & \mbox {for} & s\!\in\!(0,2.4649] \subset \left(s_1,s_2\right)       \\[1.0 ex]
\Longrightarrow & P_6(x)\!>\!0 & \mbox {for} & x\!\in\!(0,1.57] \subset \left(0,\sqrt{s_2}\right)    \\[1.0 ex]
\Longrightarrow & P_7(x)\!>\!0 & \mbox {for} & x\!\in\!(0,1.57] \subset \left(0 ,\sqrt{s_2}\right).
\end{array}
\end{equation}
where $\sqrt{s_2}=1.589\ldots >1.57\,$.

\noindent
According to the Lemmas 2.3. and 2.4. and description of the method based on (\ref{Ineq_14}) and (\ref{Ineq_17}),
the following inequalities:
$\cos y \!<\! \overline T_{k}^{\cos ,0}(y)\;(k\!=\!20)$,
$\cos y \!<\! \overline T_{k}^{\cos ,0}(y)\;(k\!=\!16)$,
$\sin y \!>\! \underline T_k^{\sin ,0}(y)\;(k\!=\!11)$ are true,
for $y\!\in\!\left(0, \sqrt{(k+3)(k+4)}\right)$. For $x\!\in\!(0,1.57]$ it is valid:
\begin{equation}
\label{Ineq_28}
f(x)
>
Q_{20}(x)
=
1-8 x^2\!-\!\overline{T}_{20}^{\cos ,0}(4 x)\!-\!8 x^2 \overline{T}_{16}^{\cos ,0}(2 x)
\!+\!P_7(x)\underline{T}_{11}^{\sin ,0}(2 x),
\end{equation}
where $Q_{20}(x)$ is the polynomial
\begin{equation}
\label{Ineq_29}
\begin{array}{rcl}
Q_{20}(x)\!\!&\!=\!&\!\! -\nfrac{16}{9280784638125} x^{10}{\Big(}262144 x^{10}-5203625 x^8+69322260 x^6     \\[2.5 ex]
\!\!&\! \!&\!\! -665557650 x^4+3412527300 x^2-5237832600{\Big )}                                              \\[0.5 ex]
\!\!&\!=\!&\!\! -\nfrac{16}{9280784638125} x^{10} Q_{10}(x),
\end{array}
\end{equation}
for $ x\!\in\!(0,1.57]$. Then, we have to determine sign of the polynomial
\begin{equation}
\label{Ineq_30}
\begin{array}{rcl}
Q_{10}(x)\!\!&\!=\!&\!\! 262144 x^{10}-5203625 x^8+69322260 x^6-665557650 x^4\\[1.0 ex]
\!\!&\! \!&\!\! +3412527300 x^2-5237832600,
\end{array}
\end{equation}
for $ x\!\in\!(0,1.57]$, which is the polynomial of $10^{th}$ degree.
By substitution $t= x^2$ we can transform the polynomial $Q_{10}(x)$ into polynomial
\begin{equation}
\label{Ineq_31}
\begin{array}{rcl}
Q_5(t)\!\!&\!=\!&\!\!262144 t^5-5203625 t^4+69322260 t^3-665557650 t^2 \\[1.0 ex]
\!\!&\! \!&\!\!+3412527300 t-5237832600,
\end{array}
\end{equation}
for $t\!\in\!(0,2.4649]$. The first derivate of the polynomial $Q_5(t)$ is the polynomial of $4^{th}$ degree
\begin{equation}
\label{Ineq_32}
Q_5^ '(t)\!=\!1312070 t^4\!-\!20814500 t^3\!+\!207966780 t^2\!-\!1331115300 t\!+\!3412527300.
\end{equation}
Using MATLAB software we can determine the real numerical factorization of the polynomial
\begin{equation}
\label{Ineq_33}
Q_5'(t)=\alpha(t^2+p_1 t+q_1)(t^2+p_2 t+q_2),
\end{equation}
where $\alpha=1310720$, $p_1=-11.655\ldots$, $q_1=34.966\ldots$, $p_2=-4.224\ldots $,
$q_2=74.457\ldots\;$. Also, holds that inequalities $p_1^2-4q_1<0$ and $p_2^2-4q_2<0$
are true. The polynomial $Q_5'(t)$ has no real roots. Let us remark that roots and
constants $p_1$, $q_1$, $p_2$, $ q_2$ can be represented in symbolic form. The polynomial
$Q_5'(t)$ is positive function for $t\!\in\!R$ therefore the polynomial $Q_5(t)$ is monotonically
increasing function for $t\!\in\!R$. Further, the function $Q_5(t)$ has real root in
$\mbox{\boldmath $a_1$}=2.464993 \ldots > 2.4649$ and $Q_5(0)<0$, so we have that
the function $Q_5(t)<0$ for $t\!\in\!(0,\mbox{\boldmath $a_1$})$ which follows
that the function $Q_{10}(x)<0$ for $ x\!\in\!(0,1.57]$.

\noindent
After all we can conclude following:
\begin{equation}
\label{Ineq_34}
\begin{array}{rccl}
                & Q_{10}(x)<0 & \mbox {for} & x\!\in\!(0,1.57] \\[1.0 ex]
\Longrightarrow & Q_{20}(x)>0 & \mbox {for} & x\!\in\!(0,1.57] \\[1.0 ex]
\Longrightarrow & f(x)>0       & \mbox {for} & x\!\in\!(0,1.57].
\end{array}
\end{equation}
Let us remark that we can easily calculate the real root $\mbox{\boldmath $a_1$}$ of the polynomial $Q_5(t)$,
and with arbitrary accuracy because $Q_5(t)$ is a strictly increasing polynomial function.
This also determines $x^*=\sqrt{a_1}=1.570029 \ldots(> 1.57)$ as the first positive root
of the polynomial $Q_{20}(x)$ defined at (\ref{Ineq_29}).

\medskip
\noindent
{\normalsize \bf (A/II)}
{\boldmath $ x\!\in\!\left(1.57,\nfrac{\pi}{2}\right)$} Let us define the function:
\begin{equation}
\label{Ineq_35}
\begin{array}{rcl}
\hat{f}(x)\!\!&\!\!=\!\!&\!\! f \left(\nfrac{\pi}{2}\!-\!x\right)\!=\!1\!-\! 8 \left(\nfrac{\pi}{2}- x \right)^{\!2}\!+\!\hat{h}_1(x)\cos4 x\!+\!\hat{h}_2(x)\cos2 x\!+\!\hat{h}_3(x)\sin2 x\\[1.5 ex]
\!\!&\!\!=\!\!&\!\! 1-8 \left(\nfrac{\pi}{2}- x \right)^{\!2}-\cos4 x+8 \left(\nfrac{\pi}{2}- x \right)^{\!2} \cos2 x \\[1.5 ex]
\!\!&\!\! \!\!&\!\!+ \left(4 \left(\nfrac{\pi}{2}- x \right)-4 \left(\nfrac{8}{45} -a \left(\nfrac{\pi}{2}- x \right) \right) \right)\sin2 x,
\end{array}
\end{equation}
where $x \!\in\!(0,c_1)$ for $c_1=\nfrac{\pi}{2}-1.57=\nfrac{\pi}{2}-\nfrac{157}{100}\left(=0.00079\ldots\right)$
and $\!\hat{h}_1(x)\!=\!-\!1\!<\!0, \hat{h}_2(x)=8\left(\nfrac{\pi}{2}- x \right)^{\!2}>0, \hat{h}_3(x)
=4\left(\nfrac{\pi}{2}- x \right)-4\left(\nfrac{8}{45}-a\left(\nfrac{\pi}{2}- x \right) \right)\!\left(\nfrac{\pi}{2}- x \right)^5$.

\medskip
\noindent
We are proving that the function $\hat{f}(x)>0$.

\noindent
Again, it is important to find sign of the polynomial $\hat{h}_3(x)$. As we see, that polynomial is the polynomial of $7^{th}$ degree or
\begin{equation}
\label{Ineq_36}
\begin{array}{rcl}
\hat{h}_3(x)=\hat{P}_7(x)\!\!&\!=\!&\!\!\!-\!\nfrac{32}{945}x^7\!+\!\nfrac{16}{135}\pi x^6\!+\!\left(\nfrac{32}{45}\!-\!\nfrac{8}{45}\pi^2\right)x^5\!+\!\left(-\nfrac{16}{9}\pi \!+\!\nfrac{4}{27}\pi^3\right)x^4 \\[1.5 ex]
\!\!&\! \!&\!\! \!+\!\left(\nfrac{16}{9}\pi^2\!-\!\nfrac{2}{27}\pi^4 \right)\!x^3\!+\!\left(-\nfrac{8}{9}\pi^3\!+\!\nfrac{1}{45}\pi^5 \right)\!x^2 \\[1.5 ex]
\!\!&\! \!&\!\!\!+\!\left(\!-\!4\!+\!\nfrac{2}{9}\pi^4\!-\!\nfrac{1}{270}\pi^6 \right)\!x+2\pi-\nfrac{1}{45} \pi^5+\nfrac{1}{3780}\pi^7.
\end{array}
\end{equation}
Using the factorization of the polynomial $\hat{P}_7(x)$ we have
\begin{equation}
\label{Ineq_37}
\begin{array}{rcl}
\hat{P}_7(x)\!\!&\!=\!&\!\! \nfrac{1}{3780}\left(-2 x+\pi \right)\Big(64 x^6-192\pi x^5+\left(240\pi^2-1344 \right) x^4                            \\[2.0 ex]
\!\!&\! \!&\!\! \!+\!\left(-\!160\pi^3\!+\!2688\pi\right)x^3\!+\!\left(60\pi^4\!-\!2016\pi^2\right)x^2\!+\!\left(-\!12\pi^5\!+\!672\pi^3\!\right)x \\[1.0 ex]
\!\!&\! \!&\!\! +\pi^6-84\pi^4+7560\Big) = \nfrac{1}{3780}\left(-2 x+\pi \right) \hat{P}_6(x),
\end{array}
\end{equation}
where
\begin{equation}
\label{Ineq_38}
\begin{array}{rcl}
\hat{P}_6(x)\!\!&\!=\!&\!\! 64 x^6-192\pi x^5+\left(240\pi^2-1344 \right) x^4+\left(-160\pi^3+2688\pi \right) x^3\\ [1.0 ex]
\!\!&\! \!&\!\! \!+\!\left(60\pi^4\!-\!2016\pi^2\right)x^2\!+\!\left(\!-\!12\pi^5\!+\!672\pi^3\right)x\!+\!\pi^6\!-\!84\pi^4\!+\!7560,
\end{array}
\end{equation}
for $ x\!\in\!(0,c_1)$. The second derivate of the polynomial $\hat{P}_6(x)$ is the polynomial of $4^{th}$ degree
\begin{equation}
\label{Ineq_39}
\begin{array}{rcl}
\hat{P}_6^{''}\!(x)\!\!&\!=\!&\!\!1920x^4\!-\!3840\pi x^3\!+\!\left(2880\pi^2\!-\!16128\right)x^2\!\\ [1.0 ex]
\!\!&\! \!&\!\! +\!\left(\!-960\pi^3\!+\!16128\pi \right)x+120\pi^4-4032\pi^2.
\end{array}
\end{equation}
Factorization of $\hat {P}_6^{''}\!(x)$ is given by
\begin{equation}
\label{Ineq_40}
\hat {P}_6^{''}\!(x)
=
24\left(20x^2\!-\!20\pi x\!+\!(5\pi^2\!-\!168)\right)\!\left(\pi\!-\!2x\right)^{2}
=
24\left(\pi\!-\!2x\right)^{2} \hat{P}_2(x),
\end{equation}
where
\begin{equation}
\label{Ineq_41}
\hat {P}_2(x)=20 x^2-20\pi x+(5\pi^2\!-168)
\end{equation}
is quadratic polynomial with two simple real roots:
\begin{equation}
\label{Ineq_42}
\hat {P}_2(x)=\alpha(x- x_1)(x- x_2),
\end{equation}
with values $ \alpha=20$, $ x_1=-1.327\ldots $, $ x_2=4.469\ldots\;$.
It holds that next inequalities are true
\begin{equation}
\label{Ineq_43}
\begin{array}{rrcl}
\hat {P}_2(x)<0 & \mbox {for} & x\!\in\!(0,c_1)\subset(x_1, x_2)            \\[1.0 ex]
\Longrightarrow\; \hat{P}_6^{''}\!(x)<0 & \mbox {for} & x\!\in\!(0,c_1)\subset(x_1, x_2).
\end{array}
\end{equation}
Therefore, for chosen interval $ x\!\in\!(0,c_1)$ the polynomial $\hat{P}_6^{{''}}(x)$ has no roots.
Since $\hat{P}_6^{{''}}(0)<0$, the polynomial $\hat{P}_6^{{''}}(x)$ is negative function for $ x\!\in\!(0,c_1)$ and $\hat{P}_6^'(x)$ is monotonically decreasing function for $ x\!\in\!(0,c_1)$.

\noindent
Furthermore, as the polynomial $\hat{P}_6^'(c_1)>0$ it follows that the polynomial $\hat{P}_6^'(x)$ is positive function for $ x\!\in\!(0,c_1)$, and the polynomial $\hat{P}_6(x)$ is monotonically increasing function for $ x\!\in\!(0,c_1)$. Because of $\hat{P}_6(0)>0$ we conclude following:
\begin{equation}
\label{Ineq_44}
\begin{array}{rccl}
& \hat{P}_6(x)>0 & \mbox {for} & x\!\in\!(0,c_1) \\[1.0 ex]
\Longrightarrow & \hat{P}_7(x)>0 & \mbox {for} & x\!\in\!(0,c_1).
\end{array}
\end{equation}
According to the Lemmas 2.3. and 2.4. and description of the method based on (14) and (17), the following inequalities: $\cos y < \overline T_{k}^{\cos ,0}(y)(k=0)$ ,  $\cos y > \underline T_{k}^{\cos ,0}(y)(k=2)$, $\sin y > \underline T_k^{\sin ,0}(y)(k=3)$ are true, for $y\!\in\!\left(0, \sqrt{(k+3)(k+4)}\right)$. For $x\!\in\!(0,c_1)$ it is valid:
\begin{equation}
\label{Ineq_45}
\begin{array}{rcl}
\hat{f}(x)\!\!&\!=\!&\!\!f\left(\nfrac{\pi}{2}- x \right)>\hat{Q}_{10}(x) = 1-8 \left(\nfrac{\pi}{2}- x \right)^{\!2}-\overline T_0^{\cos ,0}(4 x) \\[1.5 ex]
\!\!&\! \!&\!\!+8\left(\nfrac{\pi}{2}- x \right)^{\!2} \underline T_2^{\cos ,0}(2 x)+\hat{P}_7(x) \underline T_3^{\sin ,0}(2 x),\\[1.0 ex]
\end{array}
\end{equation}
where $\hat{Q}_{10}(x)$ is the polynomial
\begin{equation}
\label{Ineq_46}
\begin{array}{rcl}
\!\!\hat{Q}_{10}(x)
\!\!&\!=\!&\!\! \nfrac{128}{2835}x^{10} - \nfrac{64}{405}\pi x^9+ \left(\nfrac{32}{135}\pi^2-\nfrac{64}{63}\right) x^8+\left(-\nfrac{16}{81}\pi^3+\nfrac{352}{135}\pi\right) x^7                                                           \\[1.5 ex]
\!\!&\! \!&\!\!+\left(\nfrac{8}{81}\pi^4-\nfrac{368}{135}\pi^2+\nfrac{64}{45}\right) x^6+\left(-\nfrac{4}{135}\pi^5+\nfrac{40}{27}\pi^3-\nfrac{32}{9}\pi \right) x^5                                                                       \\[1.5 ex]
\!\!&\! \!&\!\!\!+\left(\nfrac{2}{405}\pi^6\!-\!\nfrac{4}{9}\pi^4\!+\!\nfrac{32}{9}\pi^2\!-\!\nfrac{32}{3}\right)x^4 \!+\!\left(\nfrac{\!-\!1}{2835}\pi^7\!+\!\!\nfrac{2}{27}\pi^5\!-\!\nfrac{16}{9}\pi^3\!+\!\nfrac{40}{3}\pi\!\right)x^3 \\[1.5 ex]
\!\!&\! \!&\!\!\!+\left(\!-\nfrac{1}{135}\pi^6\!+\!\nfrac{4}{9}\pi^4\!-\!4\pi^2\!-\!8\!\right)x^2+\left(\nfrac{1}{1890}\pi^7-\nfrac{2}{45}\pi^5+4\pi \right) x                                                                           \\[1.5 ex]
\!\!&\! \!&\!\! = -\nfrac{1}{5670} x\left(-2 x+\pi \right) \hat{Q}_8(x).
\end{array}
\end{equation}
Then, we have to determine sign of the polynomial
\begin{equation}
\label{Ineq_47}
\begin{array}{rcl}
\hat{Q}_8(x)\!\!&\!=\!&\!\!128x^8\!-\!384\pi x^7\!+\!\left(480\pi^2\!-\!2880\right)x^6\!+\!\left(\!-\!320\pi^3\!+\!5952\pi \right)x^5 \\ [1.0 ex]
\!\!&\! \!&\!\!\!+\!\left(120\pi^4\!-\!4752\pi^2\!+\!4032 \right)x^4\!+\!\left(\!-\!24\pi^5\!+\!1824\pi^3\!-\!8064\pi \right)x^3\\ [1.0 ex]
\!\!&\! \!&\!\!\!+\!\left(2\pi^6\!-\!348\pi^4\!+\!6048\pi^2\!-\!30240\right)x^2\!+\!\left(36\pi^5\!-\!2016\pi^3\!+\!22680\pi \right)x\\[1.0 ex]
\!\!&\!\!&\!\! -3\pi^6+252\pi^4-22680,
\end{array} \end{equation}
for $ x\!\in\!(0,c_1)$. The fourth derivate of the polynomial $\hat {Q}_8(x)$ is the polynomial of $4^{th}$ degree
\begin{equation}
\label{Ineq_48}
\begin{array}{rcl}
\hat {Q}_8^{(iv)}(x)\!\!&\!=\!&\!\! 215040 x^4-322560\pi x^3+(172800\pi^2-1036800) x^2\\[1.0 ex]
\!\!&\! \!&\!\! +(-38400\pi^3+714240\pi) x+2880\pi^4-114048\pi^2+96768.
\end{array}
\end{equation}
Using MATLAB software we can determine the real numerical factorization of the polynomial
\begin{equation}
\label{Ineq_49}
\begin{split}
\hat {Q}_8^{(iv)}(x)=\alpha(x- x_1)(x- x_2)(x- x_3)(x- x_4),
\end{split}
\end{equation}
with values $\alpha= 2.15\ldots 10^5$, $ x_1=-0.976\ldots $, $ x_2=0.674\ldots $, $ x_3=1.505\ldots $, $ x_4=3.509\ldots $
The polynomial $\hat {Q}_8^{(iv)}(x)$ has exactly four simple real roots with a symbolic radical representation and the corresponding numerical values: $ x_1$, $ x_2$, $ x_3$, $ x_4$. Therefore, the polynomial $\hat {Q}_8^{(iv)}(x)$ has no roots for $ x\!\in\!(0,c_1)$.
Since $\hat {Q}_8^{(iv)}(0)<0$, we can conclude that the polynomial $\hat {Q}_8^{(iv)}(x)$ is negative function for $ x\!\in\!(0,c_1)$, which follows that the function $\hat{Q}_8^{{''}'}(x)$ is monotonically decreasing for $ x\!\in\!(0,c_1)$. Doing the same procedure for all derivates up to $\hat{Q}_8^{'}(x)$ we have the following:
\begin{equation}
\begin{array}{lccccccr}
\hat{Q}_8^{{'''}}(c_1)\!>\!0: \!&\!\hat{Q}_8^{{'''}}(x)\!>\!0 \!&\! \mbox{for} \!&\! x\!\in\!(0, c_1) \!&\! \Longrightarrow\; \!&\! \hat{Q}_8^{{''}}(x)\nearrow \!&\! \mbox{for} \!&\! x\!\in\!(0, c_1), \\ [2.0 ex]
\hat{Q}_8^{{''}}(c_1)\!<\!0: \!&\!\hat{Q}_8^{{''}}(x)\!<\!0 \!&\! \mbox{for} \!&\! x\!\in\!(0, c_1) \!&\! \Longrightarrow\; \!&\! \hat{Q}_8^'(x) \searrow \!&\! \mbox{for} \!&\! x\!\in\!(0, c_1),       \\ [2.0 ex]
\hat{Q}_8^{'}(c_1)\!>\!0: \!&\!\hat{Q }_8^{'}(x)\!>\!0 \!&\! \mbox{for} \!&\! x\!\in\!(0,c_1) \!&\! \Longrightarrow\; \!&\! \hat{Q}_8(x) \nearrow \!&\! \mbox{for} \!&\! x\!\in\!(0, c_1).
 \end{array}
\end{equation}
After all, we conclude following:
\begin{equation}
\label{Ineq_50}
\begin{array}{rccl}
&\hat {Q}_8(x)<0 & \mbox {for} & x\!\in\!(0, c_1) \\[1.0 ex]
\Longrightarrow\; & \hat {Q}_{10}(x)>0 & \mbox {for} & x\!\in\!(0,c_1)\\[1.0 ex]
\Longrightarrow\; & \hat {f}(x)=f\left(\nfrac{\pi}{2}- x \right)>0 & \mbox {for} & x\!\in\!(0, c_1)\\[1.0 ex]
\Longrightarrow\; & f(x)>0 & \mbox {for} & x\!\in\! \left(1.57,\nfrac{\pi}{2} \right).
\end{array}
\end{equation}
Hence we proved that the function $f(x)$ is positive on interval $ x\!\in\!(0,1.57]$ we conclude that the function $f(x)$ is positive on whole interval $x\!\in\!\left(0, \nfrac{\pi}{2}\right)$.

\bigskip

{\normalsize \bf (B)}
Let us now prove the right side of inequality. If we write inequality in the following form
\begin{equation}
\label{Ineq_51}
\left(\nfrac{\sin x }{ x} \right)^{\!2}+\frac{\tan x}{ x}<2+\left(\frac{8}{45}-b_1(x) \right) x^3 \tan x \quad \mbox {for} \quad x\!\in\! \left(0,\frac{\pi}{2} \right),
\end{equation}
where \begin{equation}
\label{Ineq_52}
\begin{split}
b_1(x) = \frac{8x^2}{945}-\frac{\mbox{\boldmath $a$} x^4}{14175},\;\mbox{\boldmath $a$} = \frac{480\pi^6-40320\pi^4+3628800}{\pi^8} =17.15041 \ldots \; .
\end{split}
\end{equation}
The inequality (\ref{Ineq_51}) is equivalent to the mixed trigonometric inequality
\begin{equation}
\label{Ineq_53}
\begin{array}{rcl}
f(x)\!\!&\!=\!&\!\! 8x^2\!-\!1+h_1(x)\cos4 x+h_2(x)\cos2 x+h_3(x)\sin2 x \\[1.0 ex]
\!\!&\!=\!&\!\! 8x^2\!-\!1\!+\!\cos4 x\!+\!8x^2\cos2 x\!+\!\left(4\left(\frac{8}{45}\!-\!b_1(x)\right)x^5\!-\!4 x\right)\sin2 x>0,
\end{array}
\end{equation}
for $x\!\in\!\left(0,\frac{\pi}{2}\right)$, and $h_1(x)=1>0$, $h_2(x)=8 x^2>0$, $h_3(x)=4\left(\frac{8}{45}-b_1(x) \right) x^5-4 x$.

\noindent
Now let us consider two cases:

\medskip
\noindent
{\normalsize \bf (B/I)}
{\boldmath $ x\!\in\!(0,1.53]$}
Let us determine sign of the polynomial $h_3(x)$. As we see, that polynomial is the polynomial of $9^{th}$ degree
\begin{equation}
\label{Ineq_55}
\begin{split}
h_3(x)\!=\!P_9(x)=4\left(\nfrac{8}{45}\!-\!b_1(x)\right)x^5\!-\!4x\!=4\nfrac{\mbox{\boldmath $a$}}{14175}x^9\!-\!\nfrac{32}{945}x^7\!+\!\!\nfrac{32}{45}x^5\!-\!4x.
\end{split}
\end{equation}
Using factorization of the polynomial $P_9(x)$ we have
\begin{equation}
\label{Ineq_56}
\begin{array}{rcl}
P_9(x)\!\!&\!=\!&\!\!-\nfrac{4}{945}\nfrac{x}{\pi^8}\left(-2x+\pi \right)\!\left(\pi+2 x \right)\left (\left(8\pi^6-672\pi^4+60480 \right) x^6 \right.\\[2.5 ex]
\!\!&\! \!&\!\! \left. +\left(-168\pi^6+15120\pi^2 \right) x^4+3780\pi^4 x^2+945\pi^6\right) \\[1.0 ex]
\!\!&\!=\!&\!\! -\nfrac{4}{945} \nfrac{ x}{\pi^8} \left(-2 x+\pi \right)\!\left(\pi+2 x \right)P_6(x),
\end{array}
\end{equation}
where
\begin{equation}
\label{Ineq_57}
P_6(x)=(8\pi^6\!-\!672\pi^4\!+\!60480)x^6\!+\!(\!-\!168\pi^6\!+\!15120\pi^2)x^4\!+\!3780\pi^4x^2\!+\!945\pi^6, \end{equation}
for $x\!\in\!(0,1.53]$.
Introducing the substitution $\!s= x^2$ we can notice that the polynomial $P_6(x)$ can be transformed into polynomial of $3^{rd}$ degree
\begin{equation}
\label{Ineq_58}
\begin{split}
P_3(s)=(8\pi^6\!-\!672\pi^4\!+\!60480)s^3\!+\!(\!-\!168\pi^6\!+\!15120\pi^2)s^2\!+\!3780\pi^4s\!+\!945\pi^6,
\end{split}
\end{equation} $s\!\in\!(0,2.3409]$. Using MATLAB software we can determine the real numerical factorization of the polynomial
\begin{equation}
\label{Ineq_59}
\begin{split}
P_3(s)= \alpha(s-s_1)(s^2+ps+q),
\end{split}
\end{equation}
where $\alpha=2712.204\ldots $, $s_1=-2.221\ldots $, $p=-6.751\ldots $, $q=150.759\ldots $ whereby the inequality $p^2-4q<0$ is true. The polynomial $P_3(s)$ has exactly one real root with a symbolic radical representation and corresponding numerical value $s_1$.
Since $P_3(0)>0$ it follows that $P_3(s)>0$ for $s\!\in\!(s_1,\infty)$, so we have following conclusions:
\begin{equation}
\label{Ineq_60}
\begin{array}{rccl}
& P_3(s)>0 & \mbox {for} & s\!\in\!(0,2.3409] \subset(s_1,\infty)\\[1.0 ex]
 \Longrightarrow\; & P_6(x)>0 & \mbox {for} & x\!\in\!(0,1.53] \subset \left(0,\nfrac{\pi}{2} \right)\\[1.0 ex]
 \Longrightarrow\; & P_9(x)<0 & \mbox {for} & x\!\in\!(0,1.53] \subset \left(0,\nfrac{\pi}{2} \right).
\end{array}
\end{equation}
According to the Lemmas 2.3. and 2.4. and description of the method based on (14) and (17), the following inequalities: $\cos y > \underline T_{k}^{\cos ,0}(y)(k=22)$, $\cos y > \underline T_{k}^{\cos ,0}(y)(k=14)$, $\sin y < \overline T_k^{\sin ,0}(y)(k=13)$ are true, for $y\!\in\!\left(0, \sqrt{(k+3)(k+4)}\right)$. For $x\!\in\!(0,1.53]$ it is valid:
\begin{equation}
\label{Ineq_61}
\begin{array}{rcl}
f(x)\!\!&\!\!>\!\!&\!\!Q_{22}(x)\!=\!8x^2\!-\!1\!+\!8x^2  \underline T_{14}^{\cos ,0}(2x)\!+\! \underline T_{22}^{\cos ,0}(4x)\!+\!P_9(x) \overline T_{13}^{\sin ,0}(2x),
\end{array}
\end{equation}
where $Q_{22}(x)$ is the polynomial
\begin{equation}
\label{Ineq_62}
\begin{array}{rcl}
Q_{22}(x)\!\!\!&\!\!=\!\!&\!\!\! \left(\nfrac{-33554432}{2143861251406875}\!+\!\nfrac{1024}{5746615875\pi^2}\!-\!\nfrac{4096}{273648375\pi^4}\!+\! \nfrac{8192}{6081075\pi^8} \right)x^{22} \\[1.5 ex]
		\!\!\!&\!\! \!\!&\!\!\! \!+\! \left(\nfrac{140032}{343732764375}\!-\! \nfrac{1024}{147349125\pi^2}\!+\! \nfrac{4096}{7016625\pi^4}\!-\! \nfrac{8192}{155925\pi^8} \right) x^{20}    \\[1.5 ex]
		\!\!\!&\!\! \!\!&\!\!\!\!+\! \left(\!-\nfrac{787456}{97692469875}\!+\!\nfrac{512}{2679075\pi^2}\!-\! \nfrac{2048}{127575\pi^4}\!+\!\nfrac{4096}{2835\pi^8} \right) x^{18}           \\[1.5 ex]
		\!\!\!&\!\! \!\!&\!\!\! \!+\! \left(\nfrac{4672}{39092625}\!-\!\nfrac{1024}{297675\pi^2}\!+\!\nfrac{4096}{14175\pi^4}\!-\! \nfrac{8192}{315\pi^8} \right) x^{16}                    \\[1.5 ex]
		\!\!\!&\!\! \!\!&\!\!\!  \!+\! \left(\!-\nfrac{5888}{5108103}\!+\!\nfrac{512}{14175\pi^2}\!-\!\nfrac{2048}{675\pi^4}\!+\!\nfrac{4096}{15\pi^8} \right) x^{14}                       \\[1.5 ex]
		\!\!\!&\!\! \!\!&\!\!\!\!+\! \left(\nfrac{2752}{467775}\!-\!\nfrac{512}{2835\pi^2}\!+\!\nfrac{2048}{135\pi^4}\!-\!\nfrac{4096}{3\pi^8} \right) x^{12}                               \\[1.5 ex]
\!\!\!&\!\! \!\!&\!\!\!\!+\! \left(\!-\nfrac{128}{14175}\!+\!\nfrac{256}{945\pi^2}\!-\! \nfrac{1024}{45\pi^4}\!+\!\nfrac{2048}{\pi^8}\right) x^{10}                                         \\[1.5 ex]
		\!\!\!&\!\!=\!\!&\!\!\! \!-\nfrac{64}{2143861251406875}\nfrac{1}{\pi^8} x^{10} Q_{12}(x).
\end{array}
\end{equation}
Then, we have to determine sign of the polynomial
\begin{equation}
\label{Ineq_63}
\begin{array}{rcl}
Q_{12}(x) \!\!\!&\!\!\!=\!\!\!&\!\!\!(524288\pi^8\!-\!5969040\pi^6\!+\!501399360\pi^4\!-\!45125942400)x^{12} \\[1.25 ex]
\!\!\!&\!\!\!+\!\!\!&\!\!\!(\!-\!13646556\pi^8\!+\!232792560\pi^6\!-\!19554575040\pi^4                       \\[1.25 ex]
\!\!\!&\!\!\!+\!\!\!&\!\!\!1759911753600)x^{10}\!+\!(270011280\pi^8\!-\!6401795400\pi^6                      \\[1.25 ex]
\!\!\!&\!\!\!+\!\!\!&\!\!\!537750813600\pi^4\!-\!48397573224000)x^ {8}\!+\!(\!-\!4003360515\pi^8             \\[1.25 ex]
\!\!\!&\!\!\!+\!\!\!&\!\!\!115232317200\pi^6\!-\!9679514644800\pi^4\!+\!871156318032000)x^{6}                \\[1.25 ex]
\!\!\!&\!\!\!+\!\!\!&\!\!\!(38612227500\pi^8\!-\!1209939330600\pi^6\!+\!101634903770400\pi^4                 \\[1.25 ex]
\!\!\!&\!\!\!-\!\!\!&\!\!\!9147141339336000)x^{4} \!+\!(\!-\!197073451575\pi^8\!+\!6049696653000\pi^6        \\[1.25 ex]
\!\!\!&\!\!\!-\!\!\!&\!\!\!508174518852000\pi^4\!+\!45735706696680000)x^{2}\!\!+\!\!302484832650\pi^8        \\[1.25 ex]
\!\!\!&\!\!\!-\!\!\!&\!\!\!9074544979500\pi^6\!+\!762261778278000\pi^4\!-\!68603560045020000,
\end{array}
\end{equation}
for $ x\!\in\!(0, 1.53 ]$, which is the polynomial of $12^{th}$ degree. Introducing the substitution $s= x^2$ we can notice that the polynomial $Q_{12}(x)$ can be transformed into polynomial of $6^{th}$ degree
\begin{equation}
\label{Ineq_64}
\begin{array}{rcl}
Q_6(s)
\!\!&\!\!\!=\!\!\!&\!\!\left(524288\pi^8\!-\!5969040\pi^6\!+\!501399360\pi^4\!-\!45125942400\right)s^{6}                         \\[1.25 ex]
\!\!&\!\!\!+\!\!\!&\!\!\left(\!-13646556\pi^8\!+\!232792560\pi^6\!-\!19554575040\pi^4\!+\!1759911753600\right)s^{5}              \\[1.25 ex]
\!\!&\!\!\!+\!\!\!&\!\!\left(270011280\pi^8\!\!-\!\!6401795400\pi^6\!\!+\!\!537750813600\pi^4\!\!-\!\!48397573224000\right)s^{4} \\[1.25 ex]
\!\!&\!\!\!+\!\!\!&\!\!\left(\!-4003360515\pi^8\!+\!115232317200\pi^6\!-\!9679514644800\pi^4 \right.                             \\[1.25 ex]
\!\!&\!\!\!+\!\!\!&\!\!\left. 871156318032000\right)s^{3}+\left(38612227500\pi^8\!-\!1209939330600\pi^6 \right.                  \\[1.25 ex]
\!\!&\!\!\!+\!\!\!&\!\!\left. 101634903770400\pi^4\!-\!9147141339336000\right)s^{2}\!+\!\left(\!-\!197073451575\pi^8 \right.     \\[1.25 ex]
\!\!&\!\!\!+\!\!\!&\!\!\left. 6049696653000\pi^6\!-\!508174518852000\pi^4\!+\!45735706696680000\right)s                          \\[1.25 ex]
\!\!&\!\!\!+\!\!\!&\!\!302484832650\pi^8\!-\!9074544979500\pi^6\!+\!762261778278000\pi^4                                         \\[1.25 ex]
\!\!&\!\!\!-\!\!\!&\!\!68603560045020000.
\end{array}
\end{equation}
for $ s\!\in\!(0, 2.3409 ]$. The second derivate of the polynomial $Q_6(x)$ is the polynomial of $4^{th}$ degree
\begin{equation}
\label{Ineq_65}
\begin{array}{rcl}
Q_6^{''}(s)\!\!&\!=\!&\!\!30(524288\pi^8-5969040\pi^6+501399360\pi^4-45125942400)s^{4}   \\[1.25 ex]
\!\!&\!+\!&\!\!20(-13646556\pi^8+232792560\pi^6-19554575040\pi^4                         \\[1.25 ex]
\!\!&\!+\!&\!\!1759911753600)s^{3}+12(270011280\pi^8-6401795400\pi^6                     \\[1.25 ex]
\!\!&\!+\!&\!\!537750813600\pi^4-48397573224000)s^{2}+6(-4003360515\pi^8                 \\[1.25 ex]
\!\!&\!+\!&\!\!115232317200\pi^6-9679514644800\pi^4+871156318032000)s                    \\[1.25 ex]
\!\!&\!+\!&\!\!77224455000\pi^8-2419878661200\pi^6+203269807540800\pi^4                  \\[1.25 ex]
\!\!&\!-\!&\!\!18294282678672000.
\end{array}
\end{equation}
Using MATLAB software we can determine the real numerical factorization of the polynomial
\begin{equation}
\label{Ineq_66}
\begin{split}
Q_6^{''}(\!s)=\alpha(s - s_1)(s - s_2)(s^2+p s+q),
\end{split}
\end{equation}
with values $\alpha= 8.853\ldots 10^{10}$, $s_1=-3.45\ldots $, $s_2=5.381\ldots $, $p=-9.49\ldots $, $q=53.32\ldots $ Also, holds that inequality $p^2-4q<0$ is true. The polynomial $Q_6^{''}(s)$ has exactly two simple real roots with a symbolic radical representation and the corresponding numerical values: $s_1$, $s_2$.
Since we have that $Q_6^{''}(0)<0$ that follows $Q_6^{''}(s)<0$ for $s\!\in\!(0, 2.3409 ] \subset(s_1, s_2)$.

\noindent
Further, the function $Q_6^{'}(s)$ is monotonically decreasing function for $s\!\in\!(0,2.3409 ]$, $Q_6^{'}(1.53)>0$ and has the first positive root for $s=2.472\ldots $ which follows $Q_6^{'}(s)>0$ for $s\!\in\!(0,2.3409 ]$.

\noindent
The function $Q_6(s)$ is monotonically increasing for $s\!\in\!(0,2.3409 ]$, has the first positive root $\mbox{\boldmath $b$}=2.358\ldots$ and holds $Q_6(1.53)<0$, which follows:
\begin{equation}
\label{Ineq_67}
\begin{array}{rccl}
& Q_6(s)<0 & \mbox {for} & s\!\in\!(0,2.3409 ] \subset(0,\mbox{\boldmath $b$}) \\[1.0 ex]
\Longrightarrow\; & Q_{12}(x)>0 & \mbox {for} & x\!\in\!(0,1.53 ]              \\[1.0 ex]
\Longrightarrow\; & Q_{22}(x)>0 & \mbox {for} & x\!\in\!(0,1.53 ]              \\[1.0 ex]
\Longrightarrow\; & f(x)>0 & \mbox {for} & x\!\in\!(0, 1.53].
\end{array}
\end{equation}

\noindent
We can easily calculate the real root $\mbox{\boldmath $b$}$ of the polynomial $Q_6(s)$, and with arbitrary accuracy because of the monotonous increasing of the polynomial function. This also applies to $x^*=\sqrt{\mbox{\boldmath $b$}}=1.53579 \ldots >1.53$ (and $x^*< \nfrac{\pi}{2}$) which is the first positive root of the polynomial $Q_{22}(x)$ defined at (\ref{Ineq_62}).

\medskip
\noindent
{\normalsize \bf (B/II)}
{\boldmath $ x\!\in\! \left(1.53,\nfrac{\pi}{2} \right)$} Let us define the function
\begin{equation}
\label{Ineq_68}
\begin{array}{rcl}
\hat{f}(x)\!\!\!&\!\!\!=\!\!\!&\!\!\!f \left(\nfrac{\pi}{2}-x\right)=8\left(\nfrac{\pi}{2}-x\right)^{\!2}\!-\!1\!+\!\hat{h}_1(x)\cos4 x\!+\!\hat{h}_2(x)\cos2 x\!+\!\hat{h}_3(x)\sin2 x \\[1.5 ex]
\!\!\!&\!\!\!=\!\!\!&\!\!\!8\left(\nfrac{\pi}{2}-x\right)^{\!2}-1+\cos4 x-8\left(\nfrac{\pi}{2}-x\right)^{\!2} \cos2 x                                                                  \\[1.5 ex]
\!\!\!&\!\!\! \!\!\!&\!\!\!+\left(4\left(\nfrac{8}{45}-b_1\left(\nfrac{\pi}{2}-x\right) \right)\!\left(\nfrac{\pi}{2}- x\right)^{\!5}-4 \left(\nfrac{\pi}{2}- x \right)\right)\sin2 x,
\end{array}
\end{equation}
where $x\!\in\!(0,c_2)$ for $c_2=\nfrac{\pi}{2}-1.53=\nfrac{\pi}{2}-\nfrac{153}{100}\left(= 0.04079\ldots\right)$, and $\hat{h}_1(x)=1>0$, $\hat{h}_2(x)=8\left(\nfrac{\pi}{2}\!-\!x\right)^{\!2}>0$, $\hat{h}_3(x)=4\left(\nfrac{8}{45}\!-\!b_1 \left(\nfrac{\pi}{2}\!-\!x\right)\right)\!\left(\nfrac{\pi}{2}\!-\!x\right)^5\!-\!4\left(\nfrac{\pi}{2}\!-\!x\right)$.

\noindent
We are proving that the function $\hat{f}(x)>0$.

\noindent
Again, it is important to find sign of the polynomial $\hat{h}_3(x)$. As we see, that polynomial is the polynomial of $9^{th}$ degree
\begin{equation}
\label{Ineq_69}
\begin{array}{rcl}
\hat{h}_3(x)
\!\!&\!=\!&\!\!\hat{P}_9(x)=4 \left(\nfrac{8}{45}-b_1 \left(\nfrac{\pi}{2}- x \right)\right)\!\left(\nfrac{\pi}{2}- x \right)^{\!5}-4 \left(\nfrac{\pi}{2}- x \right)                                                         \\[1.5 ex]
\!\!&\!=\!&\!\!\nfrac {32}{45} \left(\nfrac{\pi}{2}- x \right)^{\!5}-\nfrac{32}{945} \left(\nfrac{\pi}{2}- x \right)^{\!7}+ \nfrac{4\mbox{\boldmath $a$}}{14175} \left(\nfrac{\pi}{2}- x \right)^{\!9}-4 \left(\nfrac{\pi}{2}- x \right).
\end{array}
\end{equation}
Let us determine the sign of the polynomial
\begin{equation}
\label{Ineq_70}
\begin{array}{rcl}
\hat{P}_9(x)\!\!&\!=\!&\!\!-\nfrac{1}{945\pi^8}\Big(x(\pi-x)(\pi-2x)(\pi^{12}-12\pi^{11}x+60\pi^{10}x^2-160 \pi^9 x^3 \\[1.25 ex]
\!\!&\! \!&\!\!+240\pi^8x^4\!-\!192\pi^7x^5\!+\!64\pi^6x^6\!-\!168\pi^{10}\!+\!1680\pi^9x\!-\!7056\pi^8x^2            \\[1.25 ex]
\!\!&\! \!&\!\!+16128\pi^7x^3\!-\!21504\pi^6x^4\!+\!16128\pi^5x^5\!-\!5376\pi^4x^6\!+\!30240\pi^6                     \\[1.25 ex]
\!\!&\! \!&\!\!-181440\pi^5x\!+\!665280\pi^4x^2\!-\!1451520\pi^3\!+\!1935360\pi^2x^4                                  \\[1.25 ex]
\!\!&\! \!&\!\!-1451520\pi x^5 +483840x^6)\Big)\!                                                                     \\[1.25 ex]
\!\!&\!=\!&\!\! \!-\!\nfrac{1}{945\pi^8}x(\pi\!-\!2x)(\pi\!-\!x)\hat{P}_6(x),
\end{array}
\end{equation}
for $ x\!\in\!(0,c_2)$ where
\begin{equation}
\label{Ineq_71}
\begin{array} {rcl}
\hat{P}_6(x)\!\!&\!\!=\!\!&\!\!(64\pi^6\!-\!5376\pi^4\!+\!483840)x^6\!+\!(\!-\!192\pi^7\!+\!16128\pi^5\!-\!1451520\pi)x^5 \\[1.25 ex]
\!\!&\!\! \!\!&\!\!+(240\pi^8-2150\pi^6+1935360\pi^2)x ^4+(-160\pi^9+16128\pi^7                                           \\[1.25 ex]
\!\!&\!\! \!\!&\!\!-1451520\pi^3)x ^3+(60\pi^{10}-7056\pi^8+665280\pi^4)x ^2+(-12\pi^{11}                                 \\[1.25 ex]
\!\!&\!\! \!\!&\!\!-1680\pi^9-181440\pi^5)x+\pi^{12}-168\pi^{10}+30240\pi^6.
\end{array}
\end{equation}
The second derivate of the polynomial $\hat{P}_6(x)$ is the polynomial of $4^{th}$ degree
\begin{equation}
\label{Ineq_72}
\begin{array}{rcl}
\hat{P}_6^{''}(x)\!\!&\!\!=\!\!&\!\!30(64\pi^6\!-\!5376\pi^4\!+\!483840)x^4\!+\!20(\!-\!192\pi^7\!+\!16128\pi^5\!\\ [1.0 ex]
\!\!&\!\! \!\!&\!\!-\!1451520\pi)x^3+12(240\pi^8-21504\pi^6+1935360\pi^2) x^2\\ [1.0 ex]
\!\!&\!\! \!\!&\!\!+6(-160\pi^9+16128\pi^7-1451520\pi^3) x +120\pi^{10}\\ [1.0 ex]
\!\!&\!\! \!\!&\!\!-14112\pi^8+1330560\pi^4.
\end{array}
\end{equation}
The polynomial $\hat{P}_6^{''}(x)$ has no real numerical roots for interval $x\!\in\!(0,c_2)$ whereby the function $\hat{P}_6^{''}(x)$ is positive function for $ x\!\in\!(0,c_2)$. That further means that the function $\hat{P}_6^'(x)$ is monotonically increasing function for $ x\!\in\!(0,c_2)$.

\noindent
The function $\hat{P}_6^{'}(x)$ has root for $x\!=\!\nfrac{\pi}{2}$, also holds that $\hat{P}_6^{'}(c_2)<0$, so we can conclude that $P_6^{'}(x)<0$ for $ x\!\in\!(0,c_2)$ and the function $P_6(x)$ is monotonically decreasing for $ x\!\in\!(0,c_2)$.
The function $P_6(x)$ has no roots for $ x\!\in\!(0,c_2)$ and $P_6(c_2)>0$ so we have the following:
\begin{equation}
\label{Ineq_73}
\begin{array}{rccl}
&\hat{P}_6(x)>0 & \mbox {for} & x\!\in\!(0,c_2)         \\[1.0 ex]
\Longrightarrow\; & \hat{P}_9(x)<0 & \mbox {for} & x\!\in\!(0,c_2).
\end{array}
\end{equation}
According to the Lemmas 2.3. and 2.4. and description of the method based on (14) and (17), the following inequalities: $\cos y > \underline  T_{k}^{\cos ,0}(y)(k=2)$ ,  $\cos y < \overline T_{k}^{\cos ,0}(y)(k=4)$, $\sin y < \overline T_k^{\sin ,0}(y)(k=1)$ are true, for $y\!\in\!\left(0, \sqrt{(k+3)(k+4)}\right)$. For $x\!\in\!(0,c_2)$ it is valid:
\begin{equation}
\label{Ineq_74}
\begin{array}{rcl}
\hat{f}(x)\!\!&\!\!=\!\!&\!\!f\!\left(\nfrac{\pi}{2}\!-\!x\right)\!>\hat{Q}_{10}(x)=\!8\left(\nfrac{\pi}{2}\!-\!x\right)^{\!2}\!-\!1\!+\! \underline T_2^{\cos ,0}(4x)\! \\[1.0 ex]
\!\!&\!\! \!\!&\!\! - 8\left(\nfrac{\pi}{2}\!-\!x\right)^{\!2} \overline T_4^{\cos ,0}(2x)+\hat{P}_9(x) \overline T_1^{\sin ,0}(2 x),
\end{array}
\end{equation}
where $\hat{Q}_{10}(x)$ is the polynomial
\begin{equation}
\label{Ineq_75}
\begin{array}{rcl}
\hat{Q}_{10}(x)\!\!&\!\!=\!\!&\!\!\left(-\nfrac{2048}{\pi^8}-\nfrac{256}{945\pi^2}+\nfrac{1024}{45\pi^4}\right) x^{10}+\left(\nfrac{128}{105\pi}-\nfrac{512}{5\pi^3}+\nfrac{9216}{\pi^7}\right) x^9 \\[1.75 ex]
\!\!&\!\! \!\!&\!\!+\left(-\nfrac{64}{27}+\nfrac{1024}{5\pi^2}-\nfrac{18432}{\pi^6}\right) x^8+\left(-\nfrac{3584}{15\pi}+\nfrac{21504}{\pi^5}+\nfrac{352}{135}\pi\right) x^7                       \\[1.75 ex]
\!\!&\!\! \!\!&\!\!+\left(\nfrac{1552}{9}-\nfrac{16128}{\pi^4}-\nfrac{16}{9}\pi^2\right) x^6+\left(\nfrac{8064}{\pi^3}+\nfrac{104}{135}\pi^3-\nfrac{3632}{45}\pi\right) x^5                         \\[1.75 ex]
\!\!&\!\! \!\!&\!\!+\left(16\!-\!\nfrac{28}{135}\pi^4\!-\!\nfrac{2688}{\pi^2}\!+\!\nfrac{1124}{45}\pi^2\right)x^4\!+\!\Big {(}\nfrac{2}{63}\pi^5\!+\!\nfrac{576}{\pi}\!-\!\nfrac{208}{45}\pi^3      \\[1.75 ex]
\!\!&\!\! \!\!&\!\!-16\pi \Big {)} x^3+\left(-72-\nfrac{2}{945}\pi^6+\nfrac{16}{45}\pi^4+4\pi^2 \right) x^2 \\[1.5 ex]
\!\!&\!\!=\!\!&\!\!-\nfrac{2}{945} \nfrac{x^2}{\pi^8}\hat{Q}_8(x).
\end{array}
\end{equation}
Then, we have to determine sign of the polynomial
\begin{equation}
\label{Ineq_76}
\begin{array}{rcl}
\hat{Q}_8(x)\!\!&\!\!=\!\!&\!\!(128\pi^6-10752\pi^4+967680) x^8+(-576\pi^7+48384\pi^5                              \\[1.5 ex]
\!\!&\!\! \!\!&\!\!-4354560\pi) x^7+(1120\pi^8-96768\pi^6+8709120\pi^2) x^6                                        \\[1.5 ex]
\!\!&\!\! \!\!&\!\!+(-1232\pi^9+112896\pi^7-10160640\pi^3) x^5+(840\pi^{10}                                        \\[1.5 ex]
\!\!&\!\! \!\!&\!\!-81480\pi^8+7620480\pi^4) x^4+(-364\pi^{11}+38136\pi^9                                          \\[1.5 ex]
\!\!&\!\! \!\!&\!\!-3810240\pi^5)x^3\!+\!(98\pi^{12}\!-\!11802\pi^{10}\!-\!7560\pi^8\!+\!1270080\pi^6)x^2          \\[1.5 ex]
\!\!&\!\! \!\!&\!\!+(\!-\!15\pi^{13}\!+\!2184\pi^{11}\!+\!7560\pi^9\!-\!272160\pi^7)x\!+\!\pi^{14}\!-\!168\pi^{12} \\[1.5 ex]
\!\!&\!\! \!\!&\!\!-1890\pi^{10}+34020\pi^8,
\end{array}
\end{equation}
for $ x\!\in\!(0,c_2)$.

\noindent
The fourth derivate of the polynomial $\hat{Q}_8(x)$ is the polynomial of $4^{th}$ degree
\begin{equation}
\label{Ineq_77}
\begin{array}{rcl}
\hat{Q}_8^{(iv)}(x)\!\!&\!\!=\!\!&\!\!1680(128\pi^6\!-\!10750\pi^4\!+\!967680)x^4\!+\!840(\!-\!576\pi^7\!+\!48384\pi^5 \\[1.25 ex]
\!\!&\!\! \!\!&\!\!-4354560\pi) x^3+360(1120\pi^8-96768\pi^6+8709120\pi^2) x^2                                         \\[1.25 ex]
\!\!&\!\! \!\!&\!\!+120(-1232\pi^9+112896\pi^7-10160640\pi^3) x +20160\pi^{10}                                         \\[1.25 ex]
\!\!&\!\! \!\!&\!\!-1955520\pi^8+182891520\pi^4.
\end{array}
\end{equation}
Using MATLAB software we can determine the real numerical factorization of the polynomial
\begin{equation}
\label{Ineq_78}
\hat{Q}_8^{(iv)}(x)= \alpha(x^2+p_1 x+q_1)(x^2+p_2 x+q_2),
\end{equation}
with values $\alpha= 7.29\ldots 10^7$, $p_1=-0.798\ldots $, $ q_1=1.417\ldots $, $p_2=-6.27\ldots $, $ q_2=11.111\ldots\;$.
Also, holds that inequalies $p_1^2-4q_1<0$ and $p_2^2-4q_2<0$ are true. The polynomial $\hat{Q}_8^{(iv)}(x)$ has no simple real roots but has two pairs of complex conjugate. Roots and constants $p_1$, $ q_1$, $p_2$, $ q_2$ can be represented in symbolic form.
The polynomial $\hat{Q}_8^{(iv)}(x)$ has no simple real roots for $x\!\in\!\left(0,\nfrac{\pi}{2}\right)$ and $\hat{Q}_8^{(iv)}(0)>0$. That means that $\hat{Q}_8^{(iv)}(x)>0$ for $ x\!\in\!(0, c_2) \subset \left(0,\nfrac{\pi}{2}\right)$ and the function $\hat{Q}_8^{'''}(x)$ is monotonically increasing for $ x\!\in\!(0,c_2)$.
Further, $\hat{Q}_8^{'''}(c_2)<0$ and the function $\hat{Q}_8^{'''}(x)$ has the first positive root $ x=1.00733 \ldots $ which follows that $\hat{Q}_8^{'''}(x)<0$ for $ x\!\in\!(0, c_2) \subset(0,1.00733 \ldots)$ and the function $\hat{Q}_8^{''}(x)$ is monotonically decreasing function for $ x\!\in\!(0,c_2)$.
$\hat{Q}_8^{''}(c_2)>0$ and the function $\hat{Q}_8^{''}(x)$ has the first positive root $ x=0.45455 \ldots $ which follows that $\hat{Q}_8^{''}(x)>0$ for $ x\!\in\!(0,c_2) \subset(0,0.45455 \ldots)$ and the function $\hat{Q}_8^{'}(x)$ is monotonically increasing function for $ x\!\in\!(0,c_2)$.
$\hat{Q}_8^{'}(0)>0$ and the function $\hat{Q}_8^{'}(x)$ has the first positive root $ x=1.16834 \ldots $ which follows that $\hat{Q}_8^{'}(x)>0$ for $ x\!\in\!(0,c_2) \subset(0,1.16834 \ldots)$ and the function $\hat{Q}_8(x)$ is monotonically increasing function for $ x\!\in\!(0, c_2)$. Since we have that $\hat{Q}_8(c_2)<0$
and the function $\hat{Q}_8(x)$ has the first positive root $ x=0.04383 \ldots $ we can conclude following:
\begin{equation}
\label{Ineq_79}
\begin{array}{rccl}
& \hat{Q}_8(x)<0 & \mbox {for} & x\!\in\!(0, c_2)\subset(0,0.04383 \ldots)                    \\[1.0 ex]
\Longrightarrow\; & \hat{Q}_{10}(x)>0 & \mbox {for} & x\!\in\!(0,c_2)                         \\[1.0 ex]
\Longrightarrow\; & \hat{f}(x)=f\left(\nfrac{\pi}{2}- x \right)>0 & \rm for & x\!\in\!(0,c_2) \\[1.0 ex]
\Longrightarrow\; & f(x)>0 & \mbox {for} & x\!\in\!\left(1.53,\nfrac{\pi}{2}\right).
\end{array}
\end{equation}
Hence we proved that the function $f(x)$ is positive on interval $ x\!\in\!(0,1.53]$ we conclude that the function $f(x)$ is positive on whole interval $ x\!\in\!\left(0,\nfrac{\pi}{2}\right)$.

}

\bigskip

\textbf{ The proof of Theorem 2.2.}

\bigskip

{\small
Transforming inequality (\ref{Ineq_8}) we have the following considerations.

\medskip
{\normalsize \bf (C)}
Let us prove the left side of the inequality
\begin{equation}
\label{Ineq_80}
2+\left(\nfrac{16}{\pi^4}+c(x)\right)x^3 \tan x<\left(\nfrac{\sin⁡ x}{x}\right)^{\!2}+ \nfrac{\tan x}{x},
\end{equation}
for $x\!\in\! \left(0,\nfrac{\pi}{2}\right)$. The inequality (\ref{Ineq_80}) is equivalent to the mixed trigonometric inequality
\begin{equation}
\label{Ineq_81}
\begin{array}{rcl}
f(x)\!\!&\!\!=\!\!&\!\!1-8 x^2 + h_1(x)\cos4 x+h_2(x)\cos2 x+h_3(x)\sin2 x\\ [1.0 ex]
\!\!&\!\!=\!\!&\!\! 1-8 x^2-\cos4 x-8 x^2 \cos2 x +\left(4 x-4 \left(\nfrac{16}{\pi^4}+c(x) \right) x^5 \right)\sin2 x>0,
\end{array}
\end{equation}
for $x\!\in\! \left(0,\nfrac{\pi}{2} \right) $, and $h_1(x)=-1<0$, $h_2(x)=-8 x^2<0$, $h_3(x)=4 x-4 \left(\nfrac{16}{\pi^4} +c(x) \right) x^5$.

\noindent
Now let us consider two cases:

\medskip
\noindent
{\normalsize \bf (C/I)}
{\boldmath $x\!\in\!(0,0.98]$} Let us determine sign of the polynomial $h_3(x)$.
As we see, that polynomial is the polynomial of $6^{th}$ degree
\begin{equation}
\label{Ineq_83}
\begin{array}{rcl}
h_3(x)\!\!&\!\!=\!\!&\!\!P_6(x)\!=\!4x\!-\!4\left(\nfrac{16}{\pi^4}+c(x)\right)x^5\!=\!\left(\nfrac{64}{\pi^5}\!-\!\nfrac{64}{\pi^3}\right)\!x^6+\!\left(\nfrac{\!-\!384}{\pi^4}\!+\!\nfrac{32}{\pi^2}\right)x^5\!+\!4x.
\end{array}
\end{equation}
Using factorization of the polynomial $P_6(x)$ we have
\begin{equation}
\label{Ineq_84}
\begin{array}{rcl}
P_6(x)\!\!&\!\!=\!\!&\!\!\nfrac{4 \left(-2 x+\pi \right)\!\left((8\pi^2-80)x^4+8\pi x^3+4\pi^2 x^2+2\pi^3 x+\pi^4 \right) x }{\pi^5} \\[1.5 ex]
\!\!&\!\!=\!\!&\!\!\nfrac{4 \left(-2 x+\pi \right) x }{\pi^5} P_4(x),
\end{array}
\end{equation}
where
\begin{equation}
\label{Ineq_85}
P_4(x)=(8\pi^2-80)x^4+8\pi x^3+4\pi^2 x^2+2\pi^3 x+\pi^4,
\end{equation}
for $ x\!\in\!(0,0.98]$. Using MATLAB software we can determine the real numerical factorization of the polynomial
\begin{equation}
\label{Ineq_86}
P_4(x)=\alpha(x- x_1)(x- x_2)(x^2+p x+q),
\end{equation}
with values $\alpha= -1.043\ldots $, $ x_1=-1.524\ldots $, $ x_2=25.663\ldots $, $p=0.046\ldots $, $q=2.387\ldots $\;.
Also, holds that inequality $p^2-4q<0$ is true. The polynomial $P_4(x)$ has exactly two simple real roots with a symbolic radical representation and corresponding numerical values $ x_1$, $ x_2$. Since $P_4(0)>0$ it follows that $P_4(x)>0$ for $ x\!\in\!(x_1, x_2)$, so we have following conclusion:
\begin{equation}
\label{Ineq_87}
\begin{array}{rccl}
                  & P_4(x)>0 & \; \mbox {for} \; & x\!\in\!(0,0.98 ] \subset(x_1, x_2) \\[1.0 ex]
\Longrightarrow\; & P_6(x)>0 & \; \mbox {for} \; & x\!\in\!(0,0.98].
\end{array}
\end{equation}
According to the Lemmas 2.3. and 2.4. and description of the method based on (14) and (17), the following inequalities: $\cos y < \overline T_{k}^{\cos ,0}(y)(k=12)$ ,  $\cos y < \overline T_{k}^{\cos ,0}(y)(k=8)$, $\sin y > \underline T_k^{\sin ,0}(y)(k=7)$ are true, for $y\!\in\!\left(0, \sqrt{(k+3)(k+4)}\right)$. For $x\!\in\!(0,0.98] $ it is valid:
\begin{equation}
\label{Ineq_88}
\begin{array}{rcl}
f(x)\!\!&\!\!>\!\!&\!\!Q_{13}(x)\!=\!1\!-\!8x^2\!-\! \overline T_{12}^{\cos ,0}(4x)\!-\!8 x^2 \overline T_8^{\cos ,0}(2x)\!+\!P_6(x) \underline T_7^{\sin ,0}(2x) ,
\end{array}
\end{equation}
where $Q_{13}(x)$ is the polynomial
\begin{equation}
\label{Ineq_89}
\begin{array}{rcl}
Q_{13}(x)\!\!&\!\!=\!\!&\!\! \left(-\nfrac{1024}{63\pi^5}+\nfrac{512}{315\pi^3}\right) x^{13}+\left(\nfrac{1024}{105\pi^4}-\nfrac{256}{315\pi^2}-\nfrac{16384}{467775}\right) x^{12}\\ [2.5 ex]
\!\!&\!\! \!\!&\!\!+\left(\nfrac{512}{3\pi^5}-\nfrac{256}{15\pi^3}\right) x^{11}+\left(-\nfrac{512}{5\pi^4}+\nfrac{128}{15\pi^2}+\nfrac{3376}{14175}\right) x^{10}\\ [2.5 ex]
\!\!&\!\! \!\!&\!\!+\left(-\nfrac{2560}{3\pi^5}+\nfrac{256}{3\pi^3}\right) x^{9}+\left(\nfrac{512}{\pi^4}-\nfrac{128}{3\pi^2}-\nfrac{64}{63}\right) x^{8}\\ [2.5 ex]
\!\!&\!\! \!\!&\!\!+\left(\nfrac{1280}{\pi^5}-\nfrac{128}{\pi^3}\right) x^{7}+\left(-\nfrac{768}{\pi^4}+\nfrac{64}{\pi^2}+\nfrac{64}{45}\right) x^{6}\\ [2.5 ex]
\!\!&\!\!=\!\!&\!\! -\nfrac{16}{467775}\nfrac{x^6}{\pi^5} Q_7(x),
\end{array}
\end{equation}
for $ x\!\in\!(0,0.98]$. Then, we have to determine sign of the polynomial
\begin{equation}
\label{Ineq_90}
\begin{array}{rcl}
Q_7(x)
\!&\!\!\!=\!\!\!&\!(-47520\pi^2+475200)x^7+(1024\pi^5+23760\pi^3-285120\pi)x^6           \\[1.25 ex]
\!&\!\!\! \!\!\!&\!+(498960\pi^2-4989600)x^5+(-6963\pi^5-249480\pi^3                     \\[1.25 ex]
\!&\!\!\! \!\!\!&\!+2993760\pi)x^4+(\!-\!2494800\pi^2\!+\!24948000)x^3\!+\!(29700\pi^5\! \\[1.25 ex]
\!&\!\!\! \!\!\!&\!+1247400\pi^3\!-\!14968800\pi)x^2+(3742200\pi^2-37422000)x            \\[1.25 ex]
\!&\!\!\! \!\!\!&\!-41580\pi^5-1871100\pi^3+22453200\pi,
\end{array}
\end{equation}
for $ x\!\in\!(0,0.98]$, which is the polynomial of $7^{th}$ degree. The third derivate of the polynomial $Q_7(x)$ is the polynomial of $4^{th}$ degree
\begin{equation}
\label{Ineq_91}
\begin{array}{rcl}
Q_7^{'''}(x)\!\!&\!\!\!=\!\!\!&\!\!210(-47520\pi^2+475200)x^4+120(1024\pi^5+23760\pi^3\\ [1.0 ex]
\!\!&\!\!\! \!\!\!&\!\!-285120\pi)x^3+60(498960\pi^2\!-\!4989600)x^2+24(-6963\pi^5\\ [1.0 ex]
\!\!&\!\!\! \!\!\!&\!\!-249480\pi^3+2993760\pi)x-14968800\pi^2+149688000.
\end{array}
\end{equation}
Using MATLAB software we can determine the real numerical factorization of the polynomial
\begin{equation}
\label{Ineq_92}
Q_7^{'''}(x)=\alpha(x- x_1)(x- x_2)(x- x_3)(x- x_4),
\end{equation}
with values $\alpha= 1.301\ldots 10^6$, $ x_1=-14.400\ldots $, $ x_2=-0.776\ldots $, $ x_3=0.174\ldots $, $ x_4=0.768\ldots\;$.
The polynomial $Q_7^{'''}(x)$ has exactly four simple real roots with a symbolic radical representation and the corresponding
numerical values $x_1, x_2, x_3,  x_4$. The polynomial $Q_7^{'''}(x)$ has two simple real roots on
$x\!\in\!\left(0,\nfrac{\pi}{2}\right)$ for $x= x_3$ and $ x= x_4$. Also holds that $Q_7^{'''}(0)>0$.
That means that $Q_7^{'''}(x)>0$ for $ x\!\in\!(0, x_3)\cup(x_4,\infty)$ and $Q_7^{'''}(x)<0$ for
$x\!\in\!(x_3, x_4)$ so the function $Q_7^{''}(x)$ is monotonically increasing for
$ x\!\in\!(0, x_3)\cup(x_4,\infty)$ and monotonically decreasing for $ x\!\in\!(x_3, x_4)$.
$Q_7^{''}(0)>0$, and $Q_7^{''}(0.98)>0$ and the function $Q_7^{''}(x)$ has no real roots on
$ x\!\in\!\left(0,\nfrac{\pi}{2}\right)$. That means that $Q_7^{''}(x)>0$ for
$x\!\in\!\left(0,\nfrac{\pi}{2}\right)$ so the function $Q_7^'(x)$ is monotonically increasing
for $ x\!\in\!\left(0,\nfrac{\pi}{2}\right)$. $Q_7^'(0)<0$, $Q_7^'(0.98)>0$ the function
$Q_7^'(x)$ has real root for $ x=0.30395 \ldots $ That means that $Q_7^'(x)<0$ for
$x\!\in\!(0,0.30395 \ldots)$ and $Q_7^'(x)>0$ for $ x\!\in\!(0.30395 \ldots ,\infty)$
so the function $Q_7(x)$ is monotonically decreasing for $ x\!\in\!(0,0.11545 \ldots)$
and monotonically increasing for $ x\!\in\!(0.11545 \ldots ,\infty)$.
$Q_7(0)<0$, $Q_7(0.98)<0$ and function $Q_7(x)$ has the first positive root
$ x=0.98609 \ldots\;$. That means that $Q_7(x)<0$ for $ x\!\in\!(0,0.98]$.
We can conclude following:
\begin{equation}
\label{Ineq_93}
\begin{array}{rccl}
&Q_7(x)<0 & \mbox {for} & x\!\in\!(0,0.98]                       \\[1.0 ex]
\Longrightarrow\; & Q_{13}(x)>0 & \mbox {for} & x\!\in\!(0,0.98] \\[1.0 ex]
\Longrightarrow\; & f(x)>0 & \mbox {for} & x\!\in\!(0,0.98].
\end{array}
\end{equation}
Let us notice that $x^*=0.98609\ldots$ is also the first positive root of
the approximation of the function $f(x)$, i.e. of the polynomial $Q_{13}(x)$,
defined at (\ref{Ineq_89}).

\medskip
\noindent
{\normalsize \bf (C/II)}
{\boldmath $ x\!\in\!\left(\!0.98,\nfrac{\pi}{2}\right)$} Let us define the function
\begin{equation}
\label{Ineq_94}
\begin{array}{rcl}
\hat{f}(x)\!\!&\!\!\!=\!\!\!&\!\!f\left(\nfrac{\pi}{2}- x \right)=1\!-\!8\left(\nfrac{\pi}{2}\!-\! x \right)^{\!2}\!+\!\hat{h}_1(x)\cos4 x\!+\!\hat{h}_2(x)\cos2 x\!+\!\hat{h}_3(x)\sin2 x \\[1.0 ex]
\!\!&\!\!\!=\!\!\!&\!\!1-8\left(\nfrac{\pi}{2}- x \right)^{\!2}-\cos4 x+8\left(\nfrac{\pi}{2}- x \right)^{\!2} \cos2 x                                                                     \\[1.0 ex]
\!\!&\!\!\! \!\!\!&\!\!+\left(4\left(\nfrac{\pi}{2}- x \right)-4\left(\nfrac{16}{\pi^4} +c\left(\nfrac{\pi}{2}- x \right)\right)\!\left(\nfrac{\pi}{2}- x \right)^{\!5} \right)\sin2 x,
\end{array}
\end{equation}
where $ x\!\in\!(0,c_3)$ for $c_3=\nfrac{\pi}{2}-0.98=\nfrac{\pi}{2}-\nfrac{49}{50} (= 0.59079\ldots$) and $\hat{h}_1(x)=-1<0$, $\hat{h}_2(x)=8 \left(\nfrac{\pi}{2}- x\right)^{\!2}>0$, $\hat{h}_3(x)=4 \left(\nfrac{\pi}{2}- x \right)-4 \left(\nfrac{16}{\pi^4} +c \left(\nfrac{\pi}{2}- x \right) \right)  \left(\nfrac{\pi}{2}- x\right)^{\!5}$.

\noindent
We are proving that the function $\hat{f}(x)>0$.

\noindent
It is important to find sign of the polynomial $\hat{h}_3(x)$. As we see, that polynomial is the polynomial of $6^{th}$ degree
\begin{equation}
\label{Ineq_95}
\begin{array}{rcl}
\hat{h}_3(x)\!\!&\!\!\!=\!\!\!&\!\!\hat{P}_6(x)=4\left(\nfrac{\pi}{2}- x \right)-4\left(\nfrac{16}{\pi^4} +c\left(\nfrac{\pi}{2}- x \right) \right)\!\left(\nfrac{\pi}{2}- x \right)^{\!5}  \\[1.5 ex]
\!\!&\!\!\!=\!\!\!&\!\!\left(\nfrac{640}{\pi^5}-\nfrac{64}{\pi^3} \right)x^6+\left(-\nfrac{1536}{\pi^4}+\nfrac{160}{\pi^2} \right)x^5+\left(\nfrac{1440}{\pi^3}-\nfrac{160}{\pi} \right)x^4 \\[1.5 ex]
\!\!&\!\!\! \!\!\!&\!\!+\left(-\nfrac{640}{\pi^2}+80 \right)x^3+\left(\nfrac{120}{\pi}-20\pi \right)x^2+\left(2\pi^2-4 \right)x.
\end{array}
\end{equation}
Using factorization of the polynomial $\hat{P}_6(x)$ we have
\begin{equation}
\label{Ineq_96}
\begin{array}{rcl}
\hat{P}_6(x)\!\!&\!\!\!=\!\!\!&\!\!\nfrac{1}{\pi^5}(2x(\pi-2x)(\pi^6-8\pi^5 x +24\pi^4 x^2 -32\pi^3 x^3+16\pi^2 x^4 -2\pi^4\\[1.0 ex]
\!\!&\!\!\! \!\!\!&\!\! +56\pi^3 x-208\pi^2 x^2+304\pi x^3 -160 x^4))= \nfrac{2x(\pi-2 x)}{\pi^5} \hat{P}_4(x),
\end{array}
\end{equation}
where
\begin{equation}
\label{Ineq_97}
\begin{array}{rcl}
\hat{P}_4(x)\!\!&\!\!\!=\!\!\!&\!\!\left(16\pi^2 -160 \right) x^4+ \left(304\pi-32\pi^3 \right) x^3+\left(24\pi^4 - 208\pi^2 \right) x^2\\[1.0 ex]
\!\!&\!\!\! \!\!\!&\!\!+\left(56\pi^3 -8\pi^5 \right) x +\pi^6 -2\pi^4,
\end{array}
\end{equation}
for $ x\!\in\!(0,c_3)$. Using MATLAB software we can determine the real numerical factorization of the polynomial
\begin{equation}
\label{Ineq_98}
\hat{P}_4(x)=\alpha(x- x_1)(x- x_2)(x^2+\!p x+\!q),
\end{equation}
where $\alpha=-2.086\ldots $, $ x_1=-24.092\ldots $, $ x_2=3.094\ldots $ , $p=-3.188\ldots $, $q=4.927\ldots $ whereby the inequality $p^2-4q<0$ is true. The polynomial $\hat{P}_4(x)$ has exactly two simple real roots with a symbolic radical representation and the corresponding numerical values $ x_1$, $ x_2$. Since we have that $\hat{P}_4(0)>0$ and knowing roots of the polynomial $\hat{P}_4(x)$ we have the following:
\begin{equation}
\label{Ineq_99}
\begin{array}{rccl}
            & \hat{P}_4(x)>0 & \; \mbox {for} \; & x\!\in\!(0,c_3) \subset(x_1, x_2)  \\[1.0 ex]
\Longrightarrow\; & \hat{P}_6(x) >0 & \; \mbox {for} \; & x\!\in\!(0,c_3).
\end{array}
\end{equation}
According to the Lemmas 2.3. and 2.4. and description of the method based on (14) and (17), the following inequalities: $\cos y < \overline T_{k}^{\cos ,0}(y)(k=8)$ ,  $\cos y > \underline T_{k}^{\cos ,0}(y)(k=6)$, $\sin y > \underline T_k^{\sin ,0}(y)(k=7)$ are true, for $y\!\in\!\left(0, \sqrt{(k+3)(k+4)}\right)$. For $x\!\in\!(0,c_3)$ it is valid:
\begin{equation}
\begin{array}{rcl}
\label{Ineq_100}
\hat{f}(x)\!\!&\!\!\!=\!\!\!&\!\!f\left(\nfrac{\pi}{2}\!-\!x\right)\!>\hat{Q}_{13}(x)=\!1\!-\!8\left(\nfrac{\pi}{2}\!-\!x\right)^{\!2}\!-\!\overline T_8^{\cos ,0}(4x)\\
\!\!&\!\!\! \!\!\!&\!\! +8\left(\nfrac{\pi}{2}-x\right)^{\!2} \underline T_6^{\cos ,0}(2 x)+\hat{P}_6(x) \underline T_7^{\sin ,0}(2 x),
\end{array}
\end{equation}
where $\hat{Q}_{13}(x)$ is the polynomial
\begin{equation}
\label{Ineq_101}
\begin{array}{rcl}
\!\!\!\!\!\!
\hat{Q}_{13}(x)
\!\!&\!\!\!=\!\!\!&\!\!\left(\!-\nfrac{1024}{63\pi^5}\!+\!\nfrac{512}{315\pi^3}\right)x^{13}\!+\!\left(\nfrac{4096}{105\pi^4}\!-\!\nfrac{256}{63\pi^2}\right)x^{12}\!+\!\left(\nfrac{512}{3\pi^5}\!-\!\nfrac{5632}{105\pi^3}\!+\!\nfrac{256}{63\pi}\right)x^{11} \\[2.25 ex]
\!\!&\!\!\! \!\!\!&\!\!\!+\!\left(\!-\nfrac{128}{63}\!+\!\nfrac{3712}{63\pi^2}\!-\!\nfrac{2048}{5\pi^4}\right)x^{10}\!+\!\left(\!-\nfrac{2560}{3\pi^5}\!+\!\nfrac{32}{63}\pi\!-\!\nfrac{320}{7\pi}\!+\!\nfrac{1408}{3\pi^3}\right)x^{9}                          \\[2.25 ex]
\!\!&\!\!\! \!\!\!&\!\!\!+\!\left(\nfrac{6016}{315}\!-\!\nfrac{16}{315}\pi^2\!-\!\nfrac{384}{\pi^2}\!+\!\nfrac{2048}{\pi^4}\right)x^{8}\!+\!\left(\nfrac{736}{3\pi}\!-\!\nfrac{2048}{\pi^3}\!+\!\nfrac{1280}{\pi^5}\!-\!\nfrac{208}{45}\pi \right)x^{7}          \\[2.25 ex]
\!\!&\!\!\! \!\!\!&\!\!\!+\!\left(\nfrac{3520}{3\pi^2}\!-\!\nfrac{3072}{\pi^4}\!+\!\nfrac{16}{45}\pi^2\!-\!\nfrac{4352}{45}\right)x^{6}\!+\!\left(\!-\nfrac{480}{\pi}\!+\!\nfrac{2880}{\pi^3}\!+\!\nfrac{64}{3}\pi \right)x^{5}                                  \\[2.25 ex]
\!\!&\!\!\! \!\!\!&\!\! \!+\!\left(\!-\nfrac{1280}{\pi^2}\!\!-\!\!\nfrac{4}{3}\pi^2 \!+\!\nfrac{416}{3}\right)x^{4}\!+\!\left(\nfrac{240}{\pi}\!-\!24\pi \right)x^{3} \\[2.2 ex]
\!\!&\!\!\! \!\!\!&\!\!= \!-\!\nfrac{4}{315\pi^5}x^3\hat{Q}_{10}(x).
\end{array}
\end{equation}
Then we have to determine sign of the polynomial
\begin{equation}
\label{Ineq_102}
\begin{array}{rcl}
\!\!\hat{Q}_{10}(x)
\!\!\!&\!\!\!=\!\!\!&\!\!\!\left(\!-128\pi^2\!\!+\!\!1280\right)x^{10}\!\!+\!\!\left(320\pi^3\!\!-\!\!3072\pi\right)x^{9}\!\!+\!\!\left(\!-320\pi^4\!\!+\!\!4224\pi^2\!\!-\!\!13440\right)x^{8} \\[1.25 ex]
\!\!\!&\!\!\!+\!\!\!&\!\!\!\left(160\pi^5\!-\!4640\pi^3\!\!+\!\!32256\pi\right)x^{7}\!\!+\!\!\left(\!-\!40\pi^6\!\!+\!\!3600\pi^4\!-\!36960\pi^2\!\!+\!\!67200\right)x^{6}                      \\[1.25 ex]
\!\!\!&\!\!\!+\!\!\!&\!\!\!\left(4\pi^7\!\!-\!\!1504\pi^5\!\!+\!\!30240\pi^3\!\!-\!\!161280\pi\right)x^{5}\!\!+\!\!\left(364\pi^6\!\!-\!\!19320\pi^4\!\!+\!\!161280\pi^2\right.                 \\[1.25 ex]
\!\!\!&\!\!\!-\!\!\!&\!\!\!\left.100800\right)x^{4}\!+\!\left(\!-28\pi^7\!+\!7616\pi^5\!-\!92400\pi^3+ 241920\pi\right)x^{3}\!+\!\left(\!-1680\pi^6\right.                                      \\[1.25 ex]
\!\!\!&\!\!\!+\!\!\!&\!\!\left.37800\pi^4\!-\!226800\pi^2\right)x^{2}\!+\!\left(105\pi^7- 10920\pi^5\!+\!100800\pi^3\right)x                                                                    \\[1.25 ex]
\!\!\!&\!\!\!+\!\!\!&\!\!1890\pi^6\!-\!18900\pi^4,
\end{array}
\end{equation}
for $ x\!\in\!(0,c_3)$ which is the polynomial of $10^{th}$ degree. The sixth derivate of the polynomial $\hat{Q}_{10}(x)$ is the polynomial of $4^{th}$ degree
\begin{equation}
\label{Ineq_103}
\begin{array}{rcl}
\hat{Q}_{10}^{(vi)}(x)\!\!&\!\!\!=\!\!\!&\!\!151200(-128\pi^2+1280)x^4+60480(320\pi^3-3072\pi)x^3 \\[1.0 ex]
\!\!&\!\!\! \!\!\!&\!\!\!+20160(-320\pi^4+4224\pi^2-13440)x^2+5040(160\pi^5-4640\pi^3             \\[1.0 ex]
\!\!&\!\!\! \!\!\!&\!\!\!+32256\pi)x-28800\pi^6+2592000\pi^4-26611200\pi^2+48384000.
\end{array}
\end{equation}
Using MATLAB software we can determine the real numerical factorization of the polynomial
\vspace*{-2mm}
\begin{equation}
\label{Ineq_104}
\hat{Q}_{10}^{(vi)}(x)=(x- x_1)(x- x_2)(x- x_3)(x- x_4),
\end{equation}
with values
$\alpha=2.523\ldots {10}^6$, $ x_1=-9.183\ldots $, $ x_2=-0.226\ldots $ , $ x_3=1.117\ldots $, $ x_4=1.796\ldots$\;.

\noindent
The polynomial $\hat{Q}_{10}^{(vi)}(x)$ has exactly four simple real roots with a symbolic radical representation and the corresponding numerical values: $ x_1$, $ x_2$, $ x_3$, $ x_4$.

\noindent
Since polynomial $\hat{Q}_{10}^{(vi)}(x)$ has root for $ x= x_3$ whereby the $\hat{Q}_{10}^{(vi)}(0)>0$ we have the following $\hat{Q}_{10}^{(vi)}(x)>0$ for $ x\!\in\!(0,c_3)\subset(0, x_3)$ and also the polynomial $\hat{Q}_{10}^{(v)}(x)$ is monotonically increasing function for $ x\!\in\!(0,c_3)$.

\noindent
Further, $Q_{10}^{(v)}(x)$ has the first positive root for $ x=0.16300 \ldots $ and $\hat{Q}_{10}^{(v)}(c_3)>0$ which gives us that $\hat{Q}_{10}^{(v)}(x)<0$ for $ x\!\in\!(0,0.16300 \ldots)$ and $\hat{Q}_{10}^{(v)}(x)>0$ for $ x\!\in\!(0.16300 \ldots ,c_3)$, also $\hat{Q}_{10}^{(iv)}(x)$ is monotonically decreasing function for $x\!\!\in\!\!(0,0.16300\!\ldots)$ and $\hat{Q}_{10}^{(iv)}(x)$ is monotonically increasing function for $ x\!\in\!(0.16300\ldots,c_3)$.

\noindent
$\hat{Q}_{10}^{(iv)}(x)$ has the first positive root for $ x=0.55589 \ldots $ and $\hat{Q}_{10}^{(iv)}(0)<0$ and $\hat{Q}_{10}^{(iv)}(c_3)>0$ which gives us that $\hat{Q}_{10}^{(iv)}(x)<0$ for $ x\!\in\!(0,0.55589 \ldots)$ and $\hat{Q}_{10}^{(iv)}(x)>0$ for $ x\!\in\!(0.55589 \ldots ,c_3)$, also $\hat{Q}_{10}^{'''}(x)$ is monotonically decreasing function for $ x\!\in\!(0,0.55589 \ldots)$ and monotonically increasing function for $ x\!\in\!(0.55589 \ldots ,c_3)$.

\noindent
$\!\hat{Q}_{10}^{'''}(x)$ has no root for $ x\!\in\!(0,c_3)$ and $\hat{Q}_{10}^{'''}(0)>0$ and $\hat{Q}_{10}^{'''}(0)>0$ which gives us that $\hat{Q}_{10}^{'''}(x)>0$ for $ x\!\in\!(0,c_3)$, also $\hat{Q}_{10}^{''}(x)$ is monotonically increasing function for $ x\!\in\!(0,c_3)$.

\noindent
$\hat{Q}_{10}^{''}(x)$ has the first positive root for $ x= 0.64192 \ldots $ and $\hat{Q}_{10}^{''}(c_3)<0$ which gives us that $\hat{Q}_{10}^{''}(x)<0$ for $ x\!\in\!(0,c_3) \subset(0,0.64192 \ldots)$, also $\hat{Q}_{10}〗^'(x)$ is monotonically decreasing function for $ x\!\in\!(0,c_3)$.

\noindent
$\hat{Q}_{10}^'(x)$ has no real root for $ x\!\in\!(0,c_3)$ and $\hat{Q}_{10}^'(c_3)>0$ which gives us that $\hat{Q}_{10}^'(x)>0$ for $ x\!\in\!(0,c_3)$, also $\hat{Q}_{10}(x)$ is monotonically increasing function for $ x\!\in\!(0,c_3)$.
$\hat{Q}_{10}(x)$ has real root $ x=0.66825 \ldots $ and $\hat{Q}_{10}(c_3)<0$ which gives us following
\begin{equation}
\label{Ineq_105}
\begin{array}{rccl}
& \hat{\varphi}_{10}(x)<0 & \mbox {for} & x\!\in\!(0,c_3)\subset(0,0.66825 \ldots)                      \\[0.75 ex]
\Longrightarrow\; & \hat{Q}_{13}(x)>0 & \mbox {for} & x\!\in\!(0,c_3)                                   \\[0.75 ex]
\Longrightarrow\; & \hat{\varphi}(x)=f\left(\nfrac{\pi}{2}- x \right)>0 & \mbox {for} & x\!\in\!(0,c_3) \\[0.75 ex]
\Longrightarrow\; & f(x)>0 & \mbox {for} & x\!\in\!\left(0.98,\nfrac{\pi}{2}\right).
\end{array}
\end{equation}
Hence we proved that the function $f(x)$ is positive for $ x\!\in\!(0,0.98]$, we conclude that the function $f(x)$ is positive for whole interval $ x\!\in\!\left(0,\nfrac{\pi}{2}\right)$

\medskip
{\normalsize \bf (D)}
Let us now prove the right side of the inequality
\begin{equation}
\label{Ineq_106}
\left(\nfrac{\sin x }{ x} \right)^{\!2}+\nfrac{\tan x}{ x}<2+ \left(\nfrac{16}{\pi^4}+d \left(x \right) \right) x^3 \tan x \quad \mbox{for} \quad x\!\in\! \left(0,\nfrac{\pi}{2} \right).
\end{equation}
The inequality (\ref{Ineq_106}) is equivalent to the mixed trigonometric inequality
\begin{equation}
\label{Ineq_107}
\begin{array}{rcl}
f(x)\!\!&\!\!\!=\!\!\!&\!\!8 x^2\!-\!1\!+\!h_1(x)\cos4 x\!+\!h_2(x)\cos2 x\!+\!h_3(x)\sin2 x                                             \\[1.0 ex]
\!\!&\!\!\!=\!\!\!&\!\! 8 x^2\!-\!1\!+\!\cos4 x\!+\!8 x^2 \cos2 x\!+\!\left(4\left(\nfrac{16}{\pi^4}\!+\!d(x) \right) x^5\!-\!4 x \right)\sin2 x>0,
\end{array}
\end{equation}
for $ x\!\in\!(0,\nfrac{\pi}{2})$, and $h_1(x)=1>0$, $h_2(x)=8 x^2>0$, $h_3(x)=4\left(\nfrac{16}{\pi^4}+d(x) \right) x^5-4 x$.

\noindent
Now, let us consider two cases:

\medskip
\noindent
{\normalsize \bf (D/I)}
{\boldmath $x\!\in\!(0,1.43]$} Let us determine sign of the polynomial $h_3(x)$. As we see, that polynomial is the polynomial of $7^{th}$ degree
\begin{equation}
\label{Ineq_109}
\begin{array}{rcl}
h_3(x)\!\!&\!\!\!=\!\!\!&\!\!P_7(x) = 4 \left(\nfrac{16}{\pi^4} \!+\!d \right) x^5\!-\!4 x\\[1.5 ex]
\!\!&\!\!\!=\!\!\!&\!\! 4 \left(\nfrac{16}{\pi^4}\!+\! \left(\nfrac{160}{\pi^5}\!-\!\nfrac{16}{\pi^3} \right)\!\left(\nfrac{\pi}{2}\!-\! x \right)
\!+\! \left(\nfrac{960}{\pi^6}\!-\!\nfrac{96}{\pi^4} \right)\!\left(\nfrac{\pi}{2}\!-\! x \right)^{\!2} \right) x^5\!-\!4 x\\[3.0 ex]
\!\!&\!\!\!=\!\!\!&\!\!\left(\nfrac{3840}{\pi^6}\!-\!\nfrac{384}{\pi^4} x^7\!+\! \left(\!-\!\nfrac{4480}{\pi^5}\!+\!\nfrac{448}{\pi^3} \right) x^6 \right)\!+\! \left(\nfrac{1344}{\pi^4}\!-\!\nfrac{128}{\pi^2} \right) x^5\!-\!4 x.
\end{array}
\end{equation}
Using factorization of the polynomial $P_7(x)$ we have
\begin{equation}
\label{Ineq_110}
\begin{array}{rcl}
P_7(x)\!\!&\!\!\!=\!\!\!&\!\!-4 x(-2 x+\pi)                                                                             \\[1.5 ex]
\!\!&\!\!\! \!\!\!&\!\!\nfrac{(32\pi^3 x^4-48\pi^2 x^5+\pi^5+2\pi^4 x+4\pi^3 x^2+8\pi^2 x^3-320\pi x^4+480 x^5)}{\pi^6} \\[1.5 ex]
\!\!&\!\!\!=\!\!\!&\!\!-\nfrac{4 x(-2 x+\pi)P_5(x)}{\pi^6},
\end{array}
\end{equation}
where
\begin{equation}
\label{Ineq_111}
P_5(x) =(480-48\pi^2) x^5 +(32\pi^3-320\pi) x^4+8\pi^2 x^3+4\pi^3 x^2+2\pi^4 x+\pi^5,
\end{equation}
for $ x\!\in\!(0,1.43]$. The first derivate of the polynomial $P_5(x)$ is the polynomial of $4^{th}$ degree
\begin{equation}
\label{Ineq_112}
P_5^'(x)= 5(480-48\pi^2) x^4 + 4(32\pi^3-320\pi) x^3+24\pi^2 x^2+8\pi^3 x+2\pi^4.
\end{equation}
Using MATLAB software we can determine the real numerical factorization of the polynomial
\begin{equation}
\label{Ineq_113}
P_5^'(x)=\alpha(x^2+p_1 x+q_1)(x^2+p_2 x+q_2),
\end{equation}
where $\alpha=31.294\ldots $, $p_1=1.004\ldots $, $q_1=0.647\ldots $, $p_2=-2.68\ldots $, $q_2=9.614\ldots $ whereby the inequalities $p_1^2-4q_1<0$ and $p_2^2-4q_2<0$ are true.

\noindent
The polynomial $P_5^'(x)$ has no real roots for interval $ x\!\in\!\left(-\nfrac{\pi}{2},\nfrac{\pi}{2}\right)$, $P_5^'(0)>0$ which gives that $P_5^'(x)>0$ for $ x\!\in\!\left(0,\nfrac{\pi}{2}\right)$, and it means that the function $P_5(x)$ is monotonically increasing function for $ x\!\in\!\left(0,\nfrac{\pi}{2}\right)$.
Further, the polynomial $P_5(x)$ also has no real roots for $ x\!\in\!\left(0,\nfrac{\pi}{2}\right)$, $P_5(0)>0$, which gives that $P_5(x)>0$ for $ x\!\in\!\left(0,\nfrac{\pi}{2}\right)$.

\noindent
Since the function $P_7(x)$ has real roots at $ x = 0$ and $ x =\nfrac{\pi}{2}$ we have the following conclusion
\begin{equation}
\label{Ineq_114}
\begin{array}{rccl}
& P_5(x)>0 & \mbox {for} & x\!\in\!(0,1.43]        \\[1.0 ex]
\Longrightarrow\; & P_7(x)<0 & \mbox {for} & x\!\in\!(0,1.43].
\end{array}
\end{equation}
According to the Lemmas 2.3. and 2.4. and description of the method based on (14) and (17), the following inequalities: $\cos y > \underline T_{k}^{\cos ,0}(y)(k=10)$ ,  $\sin y < \overline T_k^{\sin ,0}(y)(k=1)$ are true, for $y\!\in\!\left(0, \sqrt{(k+3)(k+4)}\right)$. For $x\!\in\!(0,1.43]$ it is valid:
\begin{equation}
\label{Ineq_115}
\begin{array}{rcl}
f(x)\!\!\!&\!\!\!>\!\!\!&\!\!\!Q_{12}(x)\!=\!8 x^2\!-\!1\!+\!\underline T_{10}^{\cos,0}(4x)\!+\!8x^2 \underline T_{10}^{\cos,0}(2x)\!+\!P_7(x) \overline T_1^{\sin,0}(2x),
\end{array}
\end{equation}
where $Q_{12}(x)$ is the polynomial
\begin{equation}
\label{Ineq_116}
\begin{array}{rcl}
Q_{12}(x)\!\!&\!\!\!=\!\!\!&\!\!-\nfrac{32}{14175}x^{12}+\left(-\nfrac{5120}{\pi^6}+\nfrac{512}{\pi^4}-\nfrac{3376}{14175}\right)x^{10}                                         \\[2.0 ex]
\!\!&\!\!\! \!\!\!&\!\!+\left(\nfrac{17920}{3\pi^5}-\nfrac{1792}{3\pi^3}\right)x^9+\left(\nfrac{7680}{\pi^6}-\nfrac{2560}{\pi^4}+\nfrac{512}{3\pi^2}+\nfrac{32}{35} \right) x^8 \\[2.0 ex]
\!\!&\!\!\! \!\!\!&\!\!+\left(-\nfrac{8960}{\pi^5}+\nfrac{896}{\pi^3} \right) x^7+\left(\nfrac{2688}{\pi^4}-\nfrac{256}{\pi^2}-\nfrac{16}{45} \right) x^6                       \\[2.0 ex]
\!\!&\!\!\!=\!\!\!&\!\! -\nfrac{16}{14175\pi^6}x^6Q_6(x),
\end{array}
\end{equation}
Then, we have to determine sign of the polynomial
\begin{equation}
\label{Ineq_117}
\begin{array}{rcl}
Q_6(x)\!\!&\!\!\!=\!\!\!&\!\! 2\pi^6 x^6\!+\!(211\pi^6\!-\!453600\pi^2\!+\!4536000) x^4\!+\!(529200\pi^3\\[1.0 ex]
\!\!&\!\!\! \!\!\!&\!\!-5292000\pi) x^3\!+\!(\!-\!810\pi^6\!-\!151200\pi^4\!+\!2268000\pi^2\!-\!6804000) x^2\\[1.0 ex]
\!\!&\!\!\! \!\!\!&\!\!+(-793800\pi^3\!+\!7938000\pi) x\!+\!315\pi^6\!+\!226800\pi^4\!-\!2381400\pi^2,
\end{array}
\end{equation}
for $ x\!\in\!(0,1.43]$. The second derivate of the polynomial $Q_6(x)$ is the polynomial of $4^{th}$ degree
\begin{equation}
\label{Ineq_118}
\begin{array}{rcl}
Q_6^{''}(x)\!\!&\!\!\!=\!\!\!&\!\! 60\pi^6 x^4+12(211\pi^6-453600\pi^2+4536000) x^2+6(529200\pi^3\\[1.0 ex]
\!\!&\!\!\! \!\!\!&\!\!-5292000\pi) x-1620\pi^6-302400\pi^4+4536000\pi^2-13608000.
\end{array}
\end{equation}
Using MATLAB software we can determine the real numerical factorization of the polynomial

\smallskip
\noindent
\begin{equation}
\label{Ineq_119}
Q_6^{''}(x)=\alpha(x^2+p_1 x+q_1)(x^2+p_2 x+q_2),
\end{equation}

\smallskip
\noindent
where $\alpha= 57683.351\ldots $, $p_1=0.413\ldots $, $q_1=54.628\ldots $, $p_2=-0.413\ldots $, $q_2=0.046\ldots $ whereby the inequalities $p_1^2-4q_1<0$ and $p_2^2-4q_2<0$ are true.

\noindent
The polynomial $Q_6^{''}(x)$ has no real roots for interval $ x\!\in\!\left(-\nfrac{\pi}{2},\nfrac{\pi}{2}\right)$, $Q_6^{''}(0)>0$ which gives that $Q_6^{''}(x)>0$ for $ x\!\in\!\left(0,\nfrac{\pi}{2}\right)$, and it means that the function $Q_6^'(x)$ is monotonically increasing function for $ x\!\in\!\left(0,\nfrac{\pi}{2}\right)$.

\noindent
Further, the polynomial $Q_6^'(x)$ also has no real roots for $ x\!\in\!\left(0,\nfrac{\pi}{2}\right)$, $Q_6^'(0)>0$, which gives that $Q_6^'(x)>0$ for $ x\!\in\!\left(0,\nfrac{\pi}{2}\right)$.

\noindent
Since the function $Q_6(x)$ has real roots at $ x = 1.436\ldots$ and $Q_6(0)=-1.108\ldots 10^6<0$ we have the following:
\begin{equation}
\label{Ineq_120}
\begin{array}{rccl}
& Q_6^'(x)>0 & \mbox {for} & x\!\in\!(0,1.43] \subset \left(0, \nfrac{\pi}{2} \right) \\[1.0 ex]
\Longrightarrow\; & Q_6(x)<0 & \mbox {for} & x\!\in\!(0,1.43]                         \\[1.0 ex]
\Longrightarrow\; & Q_{12}(x)>0 & \mbox {for} & x\!\in\!(0,1.43]                      \\[1.0 ex]
\Longrightarrow\; & f(x)>0 & \mbox {for} & x\!\in\!(0,1.43].
\end{array}
\end{equation}
Let us notice that $x^*=1.43649\ldots$ is also the first positive root of the approximation of the function $f(x)$, i.e. of the polynomial $Q_{12}(x)$, defined at (\ref{Ineq_116}).

\medskip
\noindent
{\normalsize \bf (D/II)}
{\boldmath $ x\!\in\!\left(1.43, \nfrac{\pi}{2} \right)$} Let us define the function
\begin{equation}
\label{Ineq_121}
\begin{array}{rcl}
\hat{f}\left(x\right)\!\!\!&\!\!\!\!=\!\!\!\!&\!\!\!f\left(\nfrac{\pi}{2}\!-\!x\right)\!=\!8\left(\nfrac{\pi}{2}\!-\! x\right)^{\!2}\!\!-\!\!1\!+\!\hat{h}_1 \left(x\right)\cos4x\!+\!\hat{h}_2 \left(x\right)\cos2x\!+\!\hat{h}_3 \left(x\right)\sin2x \\[2.0 ex]
\!\!&\!\!\!=\!\!\!&\!\!8\left(\nfrac{\pi}{2}\!-\! x\right)^{\!2}\!-\!1\!+\!\cos4 x\!-\!8\left(\nfrac{\pi}{2}\!-\! x\right)^{\!2} \cos2 x\\[2.0 ex]
\!\!&\!\!\! \!\!\!&\!\!\!+\!\left(4\left(\nfrac{16}{\pi^4} \!+\!d\left(\nfrac{\pi}{2}\!-\! x\right)\right)\!\left(\nfrac{\pi}{2}\!-\! x\right)^5\!-\!4\left(\nfrac{\pi}{2}\!-\! x\right)\right)\sin 2x,
\end{array}
\end{equation}
where $ x\!\in\!\left(0,c_4 \right)$ for $c_4= \nfrac{\pi}{2}\!-\!1.43=\nfrac{\pi}{2}\!-\!\nfrac{143}{100} (= 0.14\ldots$) and $\hat{h}_1 \left(x\right)=1>0$, $\hat{h}_2 \left(x\right)=8\left(\nfrac{\pi}{2}\!-\! x\right)^{\!2}\!-\!1>0$, $\hat{h}_3 \left(x\right)=4\left(\nfrac{16}{\pi^4} \!+\!d\left(\nfrac{\pi}{2}\!-\! x\right)\right)\!\left(\nfrac{\pi}{2}\!-\! x\right)^5\!-\!4\left(\nfrac{\pi}{2}\!-\! x\right)$.

\medskip
\noindent
We are proving that the function $\hat{f}(x)>0$.

\noindent
Further, it is important to find sign of the polynomial $\hat{h}_3(x)$. As we see, that polynomial is the polynomial of $7^{th}$ degree
\begin{equation}
\label{Ineq_122}
\begin{array}{rcl}
\hat{h}_3(x)\!\!&\!\!\!=\!\!\!&\!\!\hat{P}_7(x)=4 \left(\nfrac{16}{\pi^4} +d \left(\nfrac{\pi}{2}- x \right) \right)\!\left(\nfrac{\pi}{2}- x \right)^5-4 \left(\nfrac{\pi}{2}- x \right)                           \\[1.75 ex]
\!\!&\!\!\!=\!\!\!&\!\!4 \left(\nfrac{16}{\pi^4} + \left(\nfrac{160}{\pi^5}-\nfrac{16}{\pi^3} \right) x+ \left(\nfrac{960}{\pi^6}-\nfrac{96}{\pi^4} \right) x^2 \right)\!\left(\nfrac{\pi}{2}- x \right)^5-2\pi+4 x \\[1.75 ex]
\!\!&\!\!\!=\!\!\!&\!\!\left(-\nfrac{3840}{\pi^6}+\nfrac{384}{\pi^4} \right) x^7+ \left(\nfrac{8960}{\pi^5}-\nfrac{896}{\pi^3} \right) x^6+ \left(-\nfrac{8064}{\pi^4}+\nfrac{800}{\pi^2} \right) x^5               \\[1.75 ex]
\!\!&\!\!\!+\!\!\!&\!\!\left(\nfrac{3360}{\pi^3}-\nfrac{320}{\pi} \right) x^4+ \left(40-\nfrac{560}{\pi^2} \right) x^3+8\pi x^2+ \left(-2\pi^2+4 \right) x.
\end{array}
\end{equation}
Using factorization of the polynomial $\hat{P}_7(x)$ we have:
\begin{equation}
\label{Ineq_123}
\begin{array}{rcl}
\hat{P}_7(x)\!\!&\!\!\!=\!\!\!&\!\!-\nfrac{1}{\pi^6}\Big(2 x(\pi\!-\!2 x)(\pi^7\!-\!2\pi^6 x\!-\!24\pi^5 x^2\!+\!112\pi^4 x^3\!-\!176\pi^3 x^4 \\[1.25 ex]
\!\!&\!\!\!\!+\!\!\!\!&\!96\pi^2 x^5\!-\!2\pi^5\!-\!4\pi^4 x\!+\!272\pi^3 x^2\!-\!1136\pi^2 x^3\!+\!1760\pi x^4\!-\!960 x^5)\Big)              \\[1.25 ex]
\!\!&\!\!\!=\!\!\!&\!\!-\nfrac{2 x(\pi\!-\!2 x)}{\pi^6}\hat{P}_5(x),
\end{array}
\end{equation}

\vspace*{-2.0 mm}

\noindent
where

\bigskip
\noindent
\begin{equation}
\label{Ineq_124}
\begin{array}{rcl}
\hat{P}_5(x)\!\!&\!\!\!=\!\!\!&\!\!(96\pi^2-960) x^5+(1760\pi-176\pi^3) x^4+(112\pi^4-1136\pi^2) x^3 \\[1.5 ex]
\!\!&\!\!\!+\!\!\!&\!\!(272\pi^3-24\pi^5) x^2-(4\pi^4+2\pi^6) x+\pi^7-2\pi^5,
\end{array}
\end{equation}
for $ x\!\in\!(0,c_4)$. The first derivate of the polynomial $\hat{P}_5(x)$ is the polynomial of $4^{th}$ degree
\begin{equation}
\label{Ineq_125}
\begin{array}{rcl}
\hat{P}_5^'(x)\!\!&\!\!\!=\!\!\!&\!\!\!5(96\pi^2\!-\!960)x^4\!+\!4(1760\pi\!-\!176\pi^3)x^3\!+\!3(112\pi^4\!-\!1136\pi^2)x^2 \\[1.5 ex]
\!\!&\!\!\! \!\!\!&\!\!\!\!+\!2(272\pi^3\!-\!24\pi^5) x\!-\!(4\pi^4\!+\!2\pi^6).
\end{array}
\end{equation}
Using MATLAB software we can determine the real numerical factorization of the polynomial
\begin{equation}
\label{Ineq_126}
\hat{P}_5^'(x)=\alpha(x^2+p_1 x+q_1)(x^2+p_2 x+q_2),
\end{equation}
where $\alpha=-62.589\ldots$, $p_1=-0.461\ldots$, $q_1=7.871\ldots$, $p_2=-4.146\ldots$, $q_2=4.693\ldots$ whereby the inequalities $p_1^2-4q_1<0$ and $p_2^2-4q_2<0$ are true.

\noindent
The polynomial $\hat{P}_5^'(x)$ has no real roots for interval $ x\!\in\!(0,c_4)$, $P_5^'(0)<0$ which gives that $\hat{P}_5^'(x)<0$ for $ x\!\in\!(0,c_4)$, and it means that the function $\hat{P}_5(x)$ is monotonically increasing function for $ x\!\in\!(0,c_4)$. Further, the polynomial $\hat{P}_5(x)$ also has no real roots for $ x\!\in\!(0,c_4)$, $P_5(0)>0$, which gives that $\hat{P}_5(x)>0$ for $ x\!\in\!(0,c_4)$.

\noindent
Since the function $\hat{P}_7(x)$ has first positive root at $ x =\nfrac{\pi}{2}$ and $\hat{P}_7(0)=0$ we have the following:
\begin{equation}
\label{Ineq_127}
\begin{array}{rccl}
& \hat{P}_5(x)>0 & \mbox {for} & x\!\in\!(0,c_4)        \\[1.0 ex]
\Longrightarrow\; & \hat{P}_7(x)<0 & \mbox {for} & x\!\in\!(0,c_4).
\end{array}
\end{equation}
According to the Lemmas 2.3. and 2.4. and description of the method based on (14) and (17), the following inequalities: $\cos y > \underline T_{k}^{\cos ,0}(y)(k=6)$ , $\cos y < \overline T_{k}^{\cos ,0}(y)(k=4)$, $\sin y < \overline T_k^{\sin ,0}(y)(k=5)$ are true, for $y\!\in\!\left(0, \sqrt{(k+3)(k+4)}\right)$. For $x\!\in\!(0,c_4)$ it is valid:
\begin{equation}
\label{Ineq_128}
\begin{array}{rcl}
\hat{f}(x)\!\!&\!\!\!=\!\!\!&\!\!\!f\!\left(\nfrac{\pi}{2}\!-\!x\right)\!>\hat{Q}_{12}(x)
=\!8\left(\nfrac{\pi}{2}\!-\!x\right)^{\!2}\!-\!1\!+\!\underline T_6^{\cos,0}(4x)       \\[2.5 ex]
\!\!&\!\!\! \!\!\!&\!\!\!\!\!-8\left(\nfrac{\pi}{2}\!-\!x\right)^{\!2}\!\overline T_4^{\cos,0}(2x)
                             +\hat{P}_7(x) \overline T_5^{\sin,0}(2 x),
\end{array}
\end{equation}
where $\hat{Q}_{12}(x)$ is the polynomial
\begin{equation}
\label{Ineq_129}
\begin{array}{rcl}
\hat{Q}_{12}(x)\!\!&\!\!\!=\!\!\!&\!\!\!\left(\!-\!\nfrac{1024}{\pi^6}+\nfrac{512}{5\pi^4} \right) x^{12}+ \left(\nfrac{7168}{3\pi^5}\!-\!\nfrac{3584}{15\pi^3} \right) x^{11}                                                                     \\[1.75 ex]
\!\!&\!\!\!+\!\!\!&\!\!\!\left(\nfrac{5120}{\pi^6}\!-\!\nfrac{13312}{5\pi^4}+\nfrac{640}{3\pi^2} \right) x^{10}+ \left(\!-\!\nfrac{35840}{3\pi^5}+\nfrac{6272}{\pi^3\pi^3}\!-\!\nfrac{256}{3\pi} \right) x^9                                       \\[1.75 ex]
\!\!&\!\!\!+\!\!\!&\!\!\!\left(\!-\!\nfrac{7680}{\pi^6}\!+\!\nfrac{11520}{\pi^4}\!-\!\nfrac{1216}{\pi^2}\!+\!\nfrac{32}{3}\right)x^8\!+\!\left(\nfrac{17920}{\pi^5}\!-\!\nfrac{6272}{\pi^3}\!+\!\nfrac{1280}{3\pi}\!+\!\nfrac{32}{15}\pi\right)x^7 \\[1.75 ex]
\!\!&\!\!\!+\!\!\!&\!\!\!\left(\!-\!\nfrac{16128}{\pi^4}+\nfrac{7040}{3\pi^2}\!-\!\nfrac{8}{15}\pi^2\!-\!\nfrac{2848}{45} \right) x^6+ \left(\!-\!\nfrac{16}{3}\pi+\nfrac{6720}{\pi^3}\!-\!\nfrac{640}{\pi} \right) x^5                            \\[1.75 ex]
\!\!&\!\!\!+\!\!\!&\!\!\!\left(\!-\!\nfrac{1120}{\pi^2}+\nfrac{4}{3}\pi^2+\nfrac{304}{3} \right) x^4
=-  \nfrac{4}{45} \nfrac{ x^4}{\pi^6} \hat{Q}_8(x).
\end{array}
\end{equation}
Then, we have to determine sign of the polynomial
\begin{equation}
\label{Ineq_130}
\begin{array}{rcl}
\hat{Q}_8(x)\!\!&\!\!\!=\!\!\!&\!\!\!(-1152\pi^2+11520) x^8+(2688\pi^3-26880\pi) x^7+(-2400\pi^4\\[1.0 ex]
\!\!&\!\!\! \!\!\!&\!\!\!+29952\pi^2-57600) x^6+(960\pi^5-23520\pi^3+134400\pi) x^5\\[1.0 ex]
\!\!&\!\!\! \!\!\!&\!\!\!+(-120\pi^6+13680\pi^4-129600\pi^2+86400) x^4+(-24\pi^7\\[1.0 ex]
\!\!&\!\!\! \!\!\!&\!\!\!-4800\pi^5+70560\pi^3-201600\pi) x^3+(6\pi^8+712\pi^6-26400\pi^4\\[1.0 ex]
\!\!&\!\!\! \!\!\!&\!\!\!+171440\pi^2) x^2+(60\pi^7+7200\pi^5-75600\pi^3) x-15\pi^8\\[1.0 ex]
\!\!&\!\!\! \!\!\!&\!\!\!-1140\pi^6+12600\pi^4,
\end{array}
\end{equation}
for $ x\!\in\!(0,c_4)$. The fourth derivate of the polynomial $\hat{Q}_8(x)$ is the polynomial of $4^{th}$ degree
\begin{equation}
\label{Ineq_131}
\begin{split}
\hat{Q}_8^{(iv)}(x)= 1680(-1152\pi^2+11520) x^4+840(2688\pi^3-26880\pi) x^3\\
+360(-2400\pi^4+29952\pi^2-57600) x^2+120(960\pi^5-23520\pi^3\\
+134400\pi) x-2880\pi^6+328320\pi^4-3110400\pi^2+2073600.
\end{split}
\end{equation}
Using MATLAB software we can determine the real numerical factorization of the polynomial
\begin{equation}
\label{Ineq_132}
\hat{Q}_8^{(iv)}(x)=\alpha(x- x_1)(x- x_2)(x^2+p x+q),
\end{equation}

\smallskip
\noindent
where $\alpha= 2.523\ldots 10^5$, $ x_1=0.627\ldots $, $ x_2=1.89\ldots $, $p=-1.146\ldots $, $q=1.963\ldots $ whereby the inequation $p^2-4q<0$ is true.

\noindent
The polynomial $\hat{Q}_8^{(iv)}(x)$ has no real roots for interval $ x\!\in\!(0,c_4)$, $\hat{Q}_8^{(iv)}(0)>0$ which gives that $\hat{Q}_8^{(iv)}(x)>0$ for $ x\!\in\!(0,c_4)$, and it means that the function $\hat{Q}_8^{'''}(x)$ is monotonically increasing function for $ x\!\in\!(0,c_4)$.

\noindent
Further, the polynomial $\hat{Q}_8^{'''}(x)$ also has no real roots for $ x\!\in\!(0,c_4)$, $\hat{Q}_8^{'''}(0)>0$, which gives that $\hat{Q}_8^{'''}(x)>0$ for $ x\!\in\!(0,c_4)$ and means that polynomial
$\hat{Q}_8^{''}(x)$ is monotonically increasing function for $ x\!\in\!(0,c_4)$.
The polynomial $\hat{Q}_8^{''}(x)$ also has no real roots for $ x\!\in\!(0,c_4)$, $\hat{Q}_8^{''}(c_4)<0$, which gives that $\hat{Q}_8^{''}(x)<0$ for $ x\!\in\!(0,c_4)$ and means that polynomial $\hat{Q}_8^'(x)$ is monotonically decreasing function for $ x\!\in\!(0,c_4)$.
The polynomial $\hat{Q}_8^'(x)$ also has no real roots for $ x\!\in\!(0,c_4)$, $\hat{Q}_8^'(c_4)>0$, which gives that $\hat{Q}_8^'(x)>0$ for $ x\!\in\!(0,c_4)$ and means that polynomial $\hat{Q}_8(x)$ is monotonically increasing function for $ x\!\in\!(0,c_4)$.
The polynomial $\hat{Q}_8(x)$ has first positive real root at $ x = 0.38641 \ldots > c_4$, $\hat{Q}_8(c_4)<0$, which gives the following:
\begin{equation}
\label{Ineq_133}
\begin{array}{rccl}
& \hat{Q}_8(x)<0 & \; \mbox {for} \; & x\!\in\!(0,c_4)                                                 \\[1.00 ex]
\Longrightarrow\; & \hat{Q}_{12}(x)>0 & \; \mbox {for} \; & x\!\in\!(0,c_4)                            \\[0.25 ex]
\Longrightarrow\; & \hat{f}(x)=f\left(\nfrac{\pi}{2}- x\right)>0 & \; \mbox {for} \; & x\!\in\!(0,c_4) \\[0.25 ex]
\Longrightarrow\; & f(x)>0 & \; \mbox {for} \; & x\!\in\! \left(1.43, \nfrac{\pi}{2} \right).
\end{array}
\end{equation}
Hence we proved that the function $f(x)$ is positive for $ x\!\in\!(0,1.43]$, we conclude that the function $f(x)$ is positive for whole interval $ x\!\in\! \left(0, \nfrac{\pi}{2} \right)$.

}

\section{Conclusion} 
\noindent
With proving Theorem 2.1. and Theorem 2.2. is proved that is possible to extend interval defined
for inequalities given in Theorem 1.1. by (\ref{Ineq_4}) and Theorem 1.2. by (\ref{Ineq_5}).
The subject of future paper work is to determine the maximum interval for which the inequalities
given in previous theorems are true.

\bigskip
\noindent
{\bf Acknowledgements.} The second author was supported in part by Serbian Ministry of Education,
Science and Technological Development, Projects ON 174032 and III 44006. The third author
was partially supported by a Grant of the Romanian National Authority for Scientific Research,
CNCS-UEFISCDI, with the Project Number PN-II-ID-PCE-2011-3-0087.

\break

\bigskip

\bigskip


\bibstyle{plain}


\begin{thebibliography}{99}  

\bibitem{[1]}
{\sc C. Mortici}, {\em A Subtly Analysis of Wilker Inequation},
Applied Mathematics and Computation {\bf 231}, 516--520, 2014.

\bibitem{[2]}
{\sc J.$\,$B. Wilker}, {\em Problem E-3306},
American Mathematical Monthly~{\bf 96}, 1989.


\bibitem{[4]}
{\sc J.S. Sumner, A.A. Jagers, M. Vowe, J. Anglesio},
{\em Inequalities Involving Trigonometric Functions},
American Mathematical Monthly {\bf 98}~(3), 264--267, 1991.


\bibitem{[6]}
{\sc C. Mortici},
{\em The Natural Approach of Wilker-Cusa-Huygens Inequalities},
Mathematical Inequalities and Applications, Vol. {\bf 14}, 535--541, 2011.

\bibitem{[7]}
{\sc B. Male\v sevi\' c, M. Makragi\' c},
{\em A Method of Proving a Class of Inequalities of Mixed Trigonometric Polynomial Functions},
arXiv:1504.08345, 2015.

\bibitem{[8]}
{\sc B. Banjac, M. Makragi\' c, B. Male\v sevi\' c},
{\em Some~notes~on~a~method for proving inequalities by computer}, Results in Mathematics,
doi 10.1007/s00025-015-0485-8, 2015.

\bibitem{[9]}
{\sc B. Male\v sevi\' c, B. Banjac, I. Jovovi\' c},
{\em A proof of two conjectures of Chao-Ping Chen for inverse trigonometric functions},
arXiv:math/ 1508.06947, 2015.

\end{thebibliography}
\end{document}